\documentclass[a4paper,10pt,fleqn]{article}

\usepackage[english]{babel}
\usepackage[fleqn,reqno]{amsmath}
\usepackage{amsfonts}
\usepackage{amstext}
\usepackage{amsbsy}
\usepackage{amsthm}
\usepackage{amssymb}
\usepackage{epsfig}

\begin{document}

\newtheorem{theorem}{Theorem}
\newtheorem{proposition}[theorem]{Proposition}
\newtheorem{lemma}[theorem]{Lemma}
\newtheorem{definition}[theorem]{Definition}

\newenvironment{prooof}{\noindent{\bf Proof :}\hfill\\}{\qed}
\newenvironment{rem}{\noindent{\bf Remark}\hfill\\}{\qed}

\renewcommand\theequation{\thesection.\arabic{equation}}

\newcommand\nsection[1]{\section{#1}\setcounter{equation}{0}}

\newcommand{\bea}{\begin{eqnarray}}
\newcommand{\eea}{\end{eqnarray}}
\newcommand{\beas}{\begin{eqnarray*}}
\newcommand{\eeas}{\end{eqnarray*}}

\newcommand{\noi}{\noindent}
\newcommand{\non}{\nonumber}
\newcommand{\disp}{\displaystyle}
\newcommand{\rgt}{\rightarrow}

\def\wt{\widetilde}
\def\wh{\widehat}
\def\[{\left[}
\def\]{\right]}
\def\<{\langle}
\def\>{\rangle}
\def\({\left(}
\def\){\right)}
\def\a{\alpha}
\def\b{\beta}
\def\g{\gamma}
\def\de{\delta}
\def\eps{\varepsilon}
\def\z{\zeta}
\def\th{\theta}
\def\k{\kappa}
\def\l{\lambda}
\def\m{\mu}
\def\n{\nu}
\def\x{\xi}
\def\r{\rho}
\def\s{\sigma}
\def\t{\tau}
\def\u{\upsilon}
\def\p{\phi}
\def\om{\omega}
\def\Om{\Omega}
\def\d{\partial}
\def\D{\Delta}
\def\o{\otimes}
\def\N{{\mathbb N}}
\def\R{{\mathbb R}}
\def\C{{\mathbb C}}
\def\Z{{\mathbb Z}}
\def\T{{T_0\,}}
\def\tau{{T_1\,}}
\def\e{\hbox{\rm e}}
\def\r{{R_0\,}}
\def\Chid{{{\cal X}_\de}}
\def\bd{{b_\de}}
\def\bbd{{{\bf b}_\de}}
\def\br0{{K_\r}}
\def\G{{\Gamma}}
\def\tr{{\rm T}}
\def\wce{{{\wt \chi}_\eps(t)}}
\def\ck{{\chi_0(q,p)}}
\def\ckk{{\chi_0(q,p)\chi_1(x,y)}}



\title{
The radiation condition at infinity
for the high-frequency Helmholtz equation with source term:
a wave packet approach
}
\author{Fran\c cois Castella}
\date{}
\maketitle

\begin{center}
IRMAR - Universit\'e de Rennes 1

Campus Beaulieu - 35042 Rennes Cedex - France

{\tt francois.castella@univ-rennes1.fr}
\end{center}

\bigskip


\tableofcontents

\bigskip

\noi
{\bf Abstract:}
\noi
We consider the high-frequency Helmholtz equation with a given
source term, and a small absorption parameter $\a>0$.
The high-frequency (or: semi-classical) parameter is $\eps>0$.
We let $\eps$ and $\a$ go to zero simultaneously. We assume that
the zero energy is non-trapping for the underlying classical flow.
We also assume that the classical trajectories starting from the origin
satisfy a transversality condition, a generic assumption.

Under these assumptions, we prove that the solution
$u^\eps$ radiates in the outgoing direction, {\bf uniformly} in $\eps$.
In particular, the function $u^\eps$, when conveniently rescaled
at the scale $\eps$ close to the origin, is shown to converge
towards the {\bf outgoing} solution of the Helmholtz equation,
with coefficients frozen at the origin. This provides a uniform
version (in $\eps$) of the limiting absorption principle.

Writing the resolvent of the Helmholtz equation as the integral in time
of the associated semi-classical Schr\"odinger propagator, our analysis relies
on the following tools: (i) For very large times, we prove and use a uniform
version of the Egorov Theorem to estimate the time integral; (ii)
for moderate times, we prove a uniform dispersive estimate that relies
on a wave-packet approach, together with the above mentioned transversality
condition; (iii) for small times, we prove that the semi-classical
Schr\"odinger operator with variable coefficients has the same dispersive
properties as in the constant coefficients case, uniformly in $\eps$.

\medskip
\noi

\medskip

\section{Introduction}

In this article, we study the asymptotics
$\eps \rgt 0^+$ in the
following scaled Helm\-holtz equation, with unknown $w^\eps$,
\bea
\label{proto}
i \eps\, \a_\eps\, w^\eps(x)
+\frac12\D_x w^\eps(x)+n^2(\eps x) w^\eps(x)
=
S\(x\) .
\eea
In this scaling, the absorption parameter $\a_\eps>0$ is small,
i.e.
\beas
&&
\a_\eps \rgt 0^+ \; \text{  as } \;  \eps \rgt 0 .
\eeas
The limiting case $\a_\eps=0^+$ is actually allowed in our analysis. Also,
the index of refraction $n^2(\eps x)$ is almost constant,
\beas
&&
n^2(\eps x)\approx n^2(0) .
\eeas
The competition between these two effects is the key difficulty of
the present
work.

In all our analysis,
the variable $x$ belongs to $\R^d$, for some $d \geq 3$.
The index of refraction $n^2(x)$ is assumed to be
given, smooth and non-negative\footnote{
Our analysis is easily extended to the case where the refraction index
is a function that changes sign.
The only really important assumption on
the sign of $n$ is $n_\infty^2>0$, see Proposition \ref{brob}.
Otherwise,
all the arguments
given in this paper are easily adapted when $n^2(x)$ changes sign,
the analysis being actually simpler when $n^2(x)$ has the wrong
sign because contribution of terms involving $\chi_\de(H_\eps)$
vanishes in that case
(see below for the notations).
}
\bea
\label{pos}
\forall x \in \R^d , \quad  n^2(x) \geq 0,
\quad
\text{ and } \;
n^2(x)\in C^\infty(\R^d) .
\eea
It is also supposed that $n^2(x)$ goes
to a constant at infinity,
\bea
\label{dec}
&&
n^2(x)=n_\infty^2+O\(\<x\>^{-\rho}\) \; \text{ as } \; x \rgt \infty ,
\eea
for some, possibly small, exponant $\rho>0$\footnote{
Here and below we use the standard
notation $\<x\>:=(1+x^2)^{1/2}$.}.
In the language of Schr\"odinger operators, this means that the potential
$n_\infty^2-n^2(x)$ is assumed to be either short-range or long range.
Finally, the source term in (\ref{proto}) uses a function
$S(x)$ that is taken sufficiently smooth and decays fast enough at infinity.
We refer to the sequel for the very assumptions we need on the
refraction index $n^2(x)$, together with the source $S$ 
(see the statement of the main Theorem below).

Upon the $L^2$-unitary rescaling
$$
w^\eps(x)=\eps^{d/2} u^\eps(\eps x) ,
$$
the study of (\ref{proto}) is naturally linked to
the analysis of the high-frequency Helm\-holtz equation,
\bea
\label{two}
&&
i \eps \a_\eps u^\eps(x)
+\frac{\eps^2}{2}\D_x u^\eps(x)+n^2(x) u^\eps(x)
=
\frac{1}{\eps^{d/2}}
S\(\frac{x}{\eps}\) ,
\eea
where the source term $S(x/\eps)$ now plays the role
of a concentration profile at the scale $\eps$. In this picture,
the difficulty now comes from the interaction
between the oscillations induced
by the source $S(x/\eps)$, and the ones due to
the semiclassical operator $\eps^2\D/2+n^2(x)$.
We give below more complete motivations for looking at the asymptotics
in (\ref{proto}) or (\ref{two}).

The goal of this article is to prove that the solution $w^\eps$ to
(\ref{proto}) converges (in the distributional sense) to the {\bf outgoing
solution} of the natural constant coefficient Helmholtz equation,
i.e.
\bea
\label{mainst}
\non
&&
\mathop{\lim}\limits_{\eps\rgt 0}
w^\eps=w^{\rm out} \, , \; \text{ where } \; w^{\rm out} \;
\text{ is defined as the solution to }
\\
&&
\qquad
i 0^+ w^{\rm out}(x)
+\frac12\D_x w^{\rm out}(x)
+n^2(0) w^{\rm out}(x)
=
S\(x\) .
\eea
In other words,
\begin{align}
\label{outg}
\non
&
w^{\rm out}=\lim_{\de\rgt 0^+}
\(i \de+ \frac12\D_x +n^2(0)\)^{-1} S
\\
&
\qquad
=
i \int_0^{+\infty}
\exp\(i t \(\frac12\D_x+n^2(0)\)\) S \; dt
.
\end{align}
It is well-known that $w^{\rm out}$ can also be defined as
the unique solution
to $(\D_x/2+n^2(0)) w^{\rm out}=S$ that satisfies the
Sommerfeld radiation condition at infinity
\bea
\label{som}
\frac{x}{\sqrt{2} |x|} \cdot \nabla_x w^{\rm out}(x)
+
i n(0) w^{\rm out}(x)
=
O\(\frac{1}{|x|^2}\) \; \text{ as } \; |x|\rgt \infty .
\eea

The main geometric
assumptions we need on the refraction index to ensure the validity
of (\ref{mainst}) are twofolds.
First, we need that the trajectories of the Hamiltonian $\xi^2/2-n^2(x)$
at the zero energy are {\bf not trapped}. This is a standard assumption
in this context. It somehow prevents accumulation of energy in bounded
regions of space. Second, it turns out that the trajectories
that really matter in our analysis,
are those that start from the origin $x=0$,
with zero energy $\xi^2/2=n^2(0)$.
In this perspective,
we need
that these trajectories
satisfy a {\bf transversality condition}: in essence,
each such ray can self-intersect, but we require that
the self-intersection is then
``tranverse'' (see assumption {\bf (H)} page \pageref{HH}, i.e.
(\ref{93}), (\ref{94}), in section \ref{mod} below).
This second assumption
prevents accumulation of energy at the origin.

We wish to emphasize that the statement (\ref{mainst}) is not
obvious.
In particular, if the transversality assumption {\bf (H)} page \pageref{HH}
is not fullfilled, our analysis shows that  (\ref{mainst}) 
becomes false in general.
We also refer to the end of this paper for
``counterexamples''.

The central difficulty is the following.
On the one hand, the vanishing absorption parameter
$\a_\eps$ in (\ref{proto}) leads to thinking that
$w^\eps$ should satisfy the Sommerfeld radiation condition at infinity
{\bf with the variable
refraction index $n^2(\eps x)$} (see (\ref{som})). Knowing that
$\lim_{|x|\rgt\infty} n^2(\eps x)=n_\infty^2$, this roughly means
that $w^\eps$ should behave like
$
\exp(i 2^{-1/2} n_\infty |x|)/|x|
$
at infinity in $x$ (in dimension $d=3$, say).
On the other hand, the almost constant refraction index $n^2(\eps x)$
in (\ref{proto}) leads to observe that $w^\eps$ naturally goes
to a solution of the Helmholtz equation {\bf with constant refraction index
$n^2(0)$}. Hoping that we may follow the absorption coefficient $\a_\eps$
continuously along the limit $\eps\rgt 0$ in $n^2(\eps x)$, the statement
(\ref{mainst}) becomes natural, and $w^\eps$ should behave like
$
\exp(i 2^{-1/2} n(0) |x|)/|x|
$
asymptotically. But, since $n(0)\neq n_\infty$ in general,
the last two statements are contradictory ...
As we see, the strong non-local effects induced by the
Helmholtz equation make the key difficulty in following the continuous
dependence
of $w^\eps$ upon both the absorption parameter $\a_\eps\rgt 0^+$ and
on the index $n^2(\eps x)\rgt n^2(0)$.

\bigskip

Let us now give some more detailed account on our motivations
for looking at the asymptotics $\eps\rgt 0$ in (\ref{proto}).

In \cite{BCKP}, the high-frequency analysis of the Helmholtz equation
with source term is performed. More precisely, the asymptotic behaviour
as $\eps \rgt 0$ of the following equation is studied\footnote{
note that we use here a slightly different
scaling than the one used in \cite{BCKP}. This a harmless
modification that is due to mere convenience.}
\bea
\label{helme}
&&
i \eps \a_\eps u^\eps(x)
+\frac{\eps^2}{2}\D_x u^\eps(x)+n^2(x) u^\eps(x)
=
\frac{1}{\eps^{d/2}}
S\(\frac{x}{\eps}\) ,
\eea
where
the variable $x$ belongs to $\R^d$, for some $d \geq 3$,
and
the index of refraction $n^2(x)$ together with 
the concentration profile
$S(x)$ are as before
(see \cite{BCKP}).
Later,
the analysis of \cite{BCKP} was extended in \cite{CPR} to more general
oscillating/concentrating source terms. The paper \cite{CPR}
studies indeed the high-frequency analysis
$\eps\rgt 0$ in
\begin{align}
\label{helmeb}
\non
&
i \eps \a_\eps u^\eps(x)
+\frac{\eps^2}{2}\D_x u^\eps(x)+n^2(x) u^\eps(x)
=
\\
&
\qquad\qquad\qquad\qquad
\frac{1}{\eps^{q}}
\int_\Gamma
S\(\frac{x-y}{\eps}\) \;
A(y) \; \exp\( i \; \frac{\phi(x)}{\eps} \)
\; d\s(y) .
\end{align}
(See also \cite{CR} for extensions - see \cite{Fou} for the case where
$n^2$ has discontinuities).
In (\ref{helmeb}), the function $S$ 
again plays the role of a concentration profile like in (\ref{helme}),
but
the concentration occurs this time
around a smooth submanifold $\Gamma \subset \R^d$
of dimension $p$ instead of a point. On the more,
the source term here includes additional oscillations through the (smooth)
amplitude $A$ and phase $\phi$.
In these notations $d\s$ denotes
the induced euclidean surface measure on the manifold $\G$, and
the rescaling
exponant $q$ depends on the dimension of $\Gamma$
together with geometric considerations, see \cite{CPR}.

Both Helmholtz equations (\ref{helme}) and (\ref{helmeb}) modelize
the propagation of a high-frequency source wave
in a medium with scaled, variable, refraction index $n^2(x)/\eps^2$.
The scaling of the index
imposes that the waves propagating in the
medium naturally have wavelength $\eps$.
On the other hand,
the source in (\ref{helme}) as well as (\ref{helmeb})
is concentrating at the scale $\eps$,
close to the origin, or close to the surface $\Gamma$.
It thus carries oscillations at the typical wavelength $\eps$.
One may think of an antenna concentrated close to a point or to a surface, 
and emmitting
waves in the whole space. The important
phenomenon that these linear equations include precisely
lies in the {\bf resonant interaction}
between the high-frequency oscillations of the source,
and the propagative modes of the medium dictated by the index $n^2/\eps^2$.
This makes one of the key difficulties of the analysis
performed in \cite{BCKP} and \cite{CPR}.

A Wigner approach is used  in \cite{BCKP} and \cite{CPR} to
treat the high-frequency asymptotics $\eps\rgt 0$. Up to a harmless
rescaling, these papers establish that the Wigner transform
$f^\eps(x,\xi)$
of $u^\eps(x)$ satisfies, in the limit $\eps\rgt 0$, the stationnary
transport equation
\bea
\label{transp}
0^+ f(x,\xi)+\xi  \cdot \nabla_x f(x,\xi)
+\nabla_x n^2(x) \cdot \nabla_\xi f(x,\xi)
= Q(x,\xi) ,
\eea
where $f(x,\xi)=\lim f^\eps(x,\xi)$ measures
the energy carried by rays located at the point $x$ in space,
with frequency $\xi \in \R^d$.
The limiting source term $Q$ in (\ref{transp})
describes quantitatively the resonant interactions mentioned above.
In the easier case
of (\ref{helme}), one has $Q(x,\xi)=\de\(\xi^2/2-n^2(0)\) \;
\de(x) \;
|\wh S(\xi)|^2$, meaning that the asymptotic source of energy
is concentrated at the origin in $x$ (this is the factor $\de(x)$),
and it only carries
resonant frequencies $\xi$ above this point
(due to $\de\(\xi^2/2-n^2(0)\)$).
A similar but more complicated value of $Q$ is obtained in the
case of (\ref{helmeb}). In any circumstance,
equation (\ref{transp}) tells us that
the energy brought by the source $Q$ is propagated
in the whole space through the transport
operator $\xi \cdot \nabla_x + \nabla_x n^2(x) \cdot \nabla_\xi $ naturally
associated with the semi-classical operator $-\eps^2 \D_x/2 -n^2(x)$.
The term $0^+ f$ 
in (\ref{transp}) specifies a radiation condition at infinity for $f$,
that is the trace, as $\eps \rgt 0$ of
the absorption coefficient $\a_\eps >0$ in (\ref{helme}) and (\ref{helmeb}).
It gives
$f$ as the outgoing solution
$$
f(x,\xi)
=
\int_0^{+\infty} Q\(X(s,x,\xi),\Xi(s,x,\xi)\) \; ds.
$$
Here $(X(s,x,\xi),\Xi(s,x,\xi))$ is the value at time $s$ of the
characteristic curve of
$\xi \cdot \nabla_x + \nabla_x n^2(x) \cdot \nabla_\xi $
starting at point $(x,\xi)$ of phase-space (see (\ref{odeh}) below).
Obtaining the radiation condition for $f$
as the limiting effect of the absorption
coefficient
$\a_\eps$ in (\ref{helme}) 
is actually the second main difficulty of the analysis
performed in \cite{BCKP} and \cite{CPR}.

It turns out that the analysis performed in \cite{BCKP} 
relies at some point on the asymptotic behaviour
of the scaled wave function $w^\eps(x)=\eps^{d/2} u^\eps(\eps x)$
that measures the oscillation/concentration
behaviour of $u^\eps$ close to the origin.
Similarly, in \cite{CPR}
one needs to rescale $u^\eps$ around any point $y\in \Gamma$,
setting
$w^\eps_y(x):=\eps^{d/2} u^\eps(y+\eps x)$ for any such $y$.
We naturally have
\beas
&&
i \eps \a_\eps w^\eps(x)
+\frac12\D_x w^\eps(x)+n^2(\eps x) w^\eps(x)
=
S\(x\) ,
\eeas
in the case of (\ref{helme}), and a similar observation
holds true in the case of (\ref{helmeb}).
Hence
the natural rescaling leads to the analysis of the prototype equation
(\ref{proto}).
Under appropriate assumptions on $n^2(x)$ and $S(x)$,
it may be proved that $w^\eps$, solution to (\ref{proto}),
is bounded in the weighted $L^2$ space
$L^2(\<x\>^{1+\de} \; d x)$, for any $\de>0$,
uniformly in $\eps$. For a fixed value of $\eps$,
such weighted estimates are consequences of
the work by Agmon, H\"ormander,
\cite{Ag}, \cite{AH}. The fact that these bounds are uniform in $\eps$
is a consequence of the recent
(and optimal) estimates established by B. Perthame and L. Vega in
\cite{PV1}, \cite{PV2} (where the weighted $L^2$ space are replaced by
a more precise homogeneous Besov-like space).
The results in \cite{PV1} and \cite{PV2} actually
need a virial condition
of the type
$2 n^2(x)+x\cdot\nabla_x n^2(x)\geq c >0$, an inequality that {\em implies}
both 
our transversality assumption {\bf (H)} page \pageref{HH},
and the non-trapping condition, i.e. the two hypothesis made
in the present paper.
We also refer to the work by N. Burq \cite{Bu},
G\'erard and Martinez \cite{GM},
T. Jecko \cite{J}, as well as Wang and Zhang \cite{WZ},
for (not optimal) bounds in a similar spirit.
Under the weaker assumptions we make in the present paper,
a weaker bound may also be obtained as 
a consequence of our analysis.
In any case, once $w^\eps$ is seen to be bounded,
it naturally possesses a weak limit $w=\lim w^\eps$
in the appropriate space. The limit $w$ clearly
satisfies in a weak sense the equation
\bea
\label{limw}
&&
\(\frac12\D_x+n^2(0)\) w(x)=S(x).
\eea
Unfortunately, equation (\ref{limw}) does not specify $w=\lim w^\eps$
in a unique way,
and it has to be supplemented with a radiation condition at infinity.
In view of the equation (\ref{proto}) satisfied by $w^\eps$,
it has been {\bf conjectured} in \cite{BCKP} and \cite{CPR}
that $\lim w^\eps$ actually satisfies
$$
\lim w^\eps=w^{\rm out},
$$
where $w^{\rm out}$ is the outgoing solution defined before.
The present paper ans\-wers the conjecture formulated in these works. It also
gives geometric conditions for the convergence
$\lim w^\eps=w^{\rm out}$ to hold.

As a final remark, let us mention that our anaylsis is purely time-dependent.
We wish to indicate that similar results than those in the present paper
were recently and independently obtained by Wang and Zhang \cite{WZ}
using a stationary approach. Note that
their analysis requires the stronger virial condition.

\bigskip

\noi
Our main theorem is the following

\medskip

\noi
{\bf Main Theorem}

\noi
{\it
Let $w^\eps$ satisfy
$
\quad
i \eps \a_\eps w^\eps(x)
+\frac12\D_x w^\eps(x)+n^2(\eps x) w^\eps(x) = S(x),
$
\;
for some sequence $\a_\eps>0$ such that $\a_\eps \rgt 0^+$ as
$\eps \rgt 0$. 
Assume that the source term $S$ belongs to the Schwartz class ${\cal
S}(\R^d)$.
Suppose also that the index of refraction
satisfies the following set of assumptions
\begin{itemize}
\item {(smoothness, decay).}
\; 
There exists an exponent $\rho>0$,
and a positive constant $n_\infty^2>0$ such that
for any multi-index $\a \in \N^d$, there
exists a constant $C_\a>0$ with
\bea
\label{decn}
&&
\Big|
\d^\a_x
\(n^2(x)-n_\infty^2\)
\Big|
\leq
C_\a \; \<x\>^{-\rho-|\a|} .
\eea
\item {(non-trapping condition)}. \;
The trajectories associated with the Hamiltonian
$\xi^2/2-n^2(x)$ are {not trapped} at the zero energy.
In other words,
any trajectory $(X(t,x,\xi),\Xi(t,x,\xi))$ solution to
\begin{align}
\non
&
\frac{\d}{\d t} X(t,x,\xi)= \Xi(t,x,\xi) 
, 
&
X(0,x,\xi)=x ,
\\
\label{odeh}
&
\frac{\d}{\d t} \Xi(t,x,\xi)= \(\nabla_x n^2\) \(X(t,x,\xi)\) 
, 
&
\Xi(0,x,\xi)=\xi ,
\end{align}
with initial datum $(x,\xi)$ such that
$
\xi^2/2-n^2(x)=0,
$
is assumed to satisfy
$$
|X(t,x,\xi)| \rgt \infty , \quad \text{ as } \; |t| \rgt \infty .
$$
\item
{(tranversality condition)}. \;
The tranvsersality condition {\bf (H)} page \pageref{HH}
(see also (\ref{93}) and (\ref{94}))
on the trajectories
starting from the origin $x=0$, with zero energy $\xi^2/2=n^2(0)$,
is satisfied.
\end{itemize}

\noi
Then,
we do have the following convergence, weakly,
when tested against any function $\phi \in {\cal S}(\R^d)$,
$$
w^\eps \rgt w^{\rm out}.
$$
}

\noi
{\bf First remark}

\noi
Still referring to
{\bf (H)} page \pageref{HH},
or (\ref{93}), (\ref{94})) for the precise statements,
we readily indicate that
the transversality assumption {\bf (H)} essentially requires
that the set 
$$
\{(\eta,\xi,t) \in \R^{2d} \times ]0,\infty[ \text{ s.t. }
X(t,0,\xi)=0, \; \Xi(t,0,\xi)=\eta, \; \xi^2/2=n^2(0)\}
$$
is a smooth submanifold of $\R^{2d+1}$,
having a codimension $>d+2$, a generic asssumption.
In other words, zero energy trajectories issued from the origin
and passing several times through the origin $x=0$ should be ``rare''.
\qed
\medskip

\noi
{\bf Second remark}

\noi
As we already mentionned, it is easily proved 
that the virial condition
$2 n^2(x)+x\cdot\nabla_xn^2(x)\geq c>0$
implies both the non-trapping and the transversality
conditions. This observation relies on the identities
$\d_t\(X(t,x,\xi)^2/2\)=X(t,x,\xi) \cdot \Xi(t,x,\xi)$ and
$\d_t\(X(t,x,\xi) \cdot \Xi(t,x,\xi)\)=
\[2 n^2(x)+x \cdot \nabla_x n^2(x)\]|_{x=X(t,x,\xi)}\geq c>0$,
where $(X(t,x,\xi),\Xi(t,x,\xi)$ is any trajectory with zero
energy (see section \ref{brbr} for computations in this spirit).

In fact, the virial condition implies even more, namely that
trajectories issued from the origin with zero energy
{\em never come back to the origin}. In other words, the set involved
in assumption {\bf (H)} page \pageref{HH} is simply {\em void},
and {\bf (H)} is trivially true under the virial condition.
As the reader may easily check, such a situation allows to
considerably simplify the proof we give here: the tools developped
in sections \ref{st}, \ref{2z}, \ref{kappa}, \ref{brbr} are actually enough
to make the complete analysis, and one does not need to
go into the detailed computations of section \ref{mod}
in that case.

Last, the above Theorem asserts the convergence of $w^\eps$:
note in passing that even
the weak boundedness of $w^\eps$ under the sole above assumptions
(i.e. without the virial condition)
is not a known result. 
\qed
\medskip

The above theorem is not only a local convergence result,
valid for test functions $\phi \in {\cal S}$.
Indeed, by density of smooth functions
in weighted $L^2$ spaces, it readily implies the following immediate corollary.
It states that, provided $w^\eps$ is bounded in the natural
weighted $L^2$ space,
the convergence also holds weakly in this space. In other words,
the convergence also holds globally.

\medskip

\noi
{\bf Immediate corollary}

\noi
{\it
With the notations of the main Theorem,
assume that the source term $S$ above satisfies the weaker decay property
\bea
&&
\label{sinb}
\|S\|_B:=\sum_{j\in \Z}
2^{j/2}
\|S\|_{L^2(C_j)} < \infty ,
\eea
where $C_j$ denotes the annulus $\{ 2^j \leq |x| \leq 2^{j+1}\}$
in $\R^d$.
Suppose also that the index of refraction
satisfies the smoothness condition of the main Theorem, with the
non-trapping and transversality assumptions replaced by the stronger
\bea
\label{pv}
&&
\bullet
\text{ (virial-like condition) } \qquad
2 \sum_{j\in\Z} \; 
\sup_{x\in C_j}
\frac{\(x \cdot \nabla n^2(x)\)_{-}}{n^2(x)}
< 1 .
\eea
Then,
we do have the convergence
$
w^\eps \rgt w^{\rm out},
$
weakly, when tested against any function $\phi$ such that
$\|\phi\|_B<\infty$,
}

\medskip

Under the simpler virial condition
$2 n^2(x)+x \cdot n^2(x)\geq c>0$,
a similar result holds with the space $B$ replaced by the more usual weighted
space  $L^2\(\<x\>^{1+\de} dx\)$ ($\de>0$ arbitrary).
Here, we give a version where
the decay (\ref{sinb}) assumed on the source $S$
is the optimal one, and the above weak convergence holds
in the optimal space.

It is well known that the resolvent of the Helmholtz operator maps
the weighted $L^2$ space $L^2\(\<x\>^{1+\de} dx\)$
to $L^2\(\<x\>^{-1-\de} dx\)$ for any $\de>0$
(\cite{Ag}, \cite{J}, \cite{GM}).
Agmon and H\"ormander \cite{AH} gave an optimal
version in the constant coefficients case:
the resolvent of the Helmholtz operator sends
the weighted $L^2$ space $B$ defined in (\ref{sinb}) to
the dual weighted space $B^*$ defined by
\begin{align}
\label{b*}
&&
\|u\|_{B^*}
:=
\sup_{j \in \Z}
2^{-j/2} \|u\|_{L^2(C_j)} .
\end{align}
For non-constant coefficients,
that are non-compact perturbations of constants,
Perthame and Vega in \cite{PV1} and \cite{PV2} established
the optimal estimate in $B$-$B^*$ under assumption
(\ref{pv}).
In our perspective, the assumption (\ref{pv}) is of
technical nature, and it may be replaced by {\bf any} assumption ensuring that
the solution $w^\eps$ to (\ref{proto}) satisfies the uniform bound
\begin{align}
\label{bbb}
&&
\|w^\eps\|_{B^*}
\leq
C_{d,n^2}  \;
\|S\|_B ,
\end{align}
for some universal constant $C_{d,n^2}$ that only depends on the dimension
$d\geq 3$ and the index $n^2$.

\medskip

\noi
{\bf Proof of the immediate Corollary}

\noi
Under the virial-like assumption (\ref{pv}), 
it has been established in \cite{PV1} that estimate
(\ref{bbb}) holds true. Hence, by density of the Schwartz class
in the space $B$, one readily reduces the problem to the case
when the source $S$ and the test function $\phi$ belong to ${\cal S}(\R^d)$.
The Main Theorem now allows to conclude.
\qed

\medskip

Needless to say, the central assumptions needed for the theorem are
the non-trapping condition together with the transversality condition.
Comments are given below on the very meaning of the transversality condition
{\bf (H)} page \pageref{HH} (i.e. (\ref{93}), (\ref{94})), to which we refer.

To state the result very briefly, the heart of our proof lies in proving that
under the above assumptions, the propagator
$\exp\(i \eps^{-1} t \, \(-\eps^2 \D_x/2-n^2(x)\)\)$,
or its rescaled value $\exp\(i t \(- \D_x/2-n^2(\eps x)\)\)$,
satisfy ``similar'' dispersive properties
as the free Schr\"odinger operator $\exp\(i t \(- \D_x/2-n^2(0)\)\)$,
{\em uniformly in~$\eps$}. This in turn is proved upon distinguishing between
small times, moderate times, and very large times,
each case leading to the use of different arguments and techniques.

\medskip

The remainder part of this paper is devoted to the proof of the main
Theorem. The proof being long and using many different tools,
we first draw in section \ref{redu} an outline of the proof,
giving the main ideas and tools. We also define the relevant
mathematical objects to be used throughout the paper.
The proof itself is performed in the next sections
\ref{st} to \ref{ccl}. Examples and counterexamples to the Theorem
are also proposed in the last section \ref{expl}.

The main intermediate results are
proposition \ref{ouane},
proposition \ref{2},
proposition \ref{wawang},
together with the more difficult proposition \ref{brob}
(that needs an Egorov Theorem for large times stated
in Lemma \ref{bourob}).
The key (and most difficult) result is  proposition \ref{dur}.
The latter uses the tranversality condition mentioned before.


\section{Preliminary Analysis: outline of the proof of the Main Theorem}
\label{redu}

\setcounter{equation}{0}

\subsection{Outline of the proof} 


Let $w^\eps$ be the solution to
$
i\eps \a_\eps w^\eps
+\frac12\D w^\eps
+n^2(\eps x) w^\eps
=
S\(x\),
$
with $S\in {\cal S}(\R^d)$.
According to the statement of our main Theorem,
we wish to study the asymptotic behaviour of $w^\eps$ as $\eps \rgt
0$, in a weak sense. Taking a test function
$\phi(x)\in {\cal S}(\R^d)$, and defining the duality product
$$
\<w^\eps,\phi\>
:=
\int_{\R^d}
w^\eps(x) \phi(x) \, dx,
$$
we want to prove the convergence
$$
\<w^\eps,\phi\>\rgt \<w^{\rm out},\phi\> \; \text{ as } \; \eps\rgt 0
.
$$
where the outgoing solution of the (constant coefficient)
Helmholtz equation $w^{\rm out}$ is defined in (\ref{mainst}), (\ref{outg})
before.

\medskip

\noi
{\bf First step: preliminary reduction - the time dependent approach}

In order to prove the weak convergence
$\<w^\eps,\phi\>\rgt\<w,\phi\>$, 
we define
the rescaled function
\bea
\label{wtou}
&&
u^\eps(x)=\frac{1}{\eps^{d/2}}
w^\eps\(\frac{x}{\eps}\)
.
\eea
It satisfies
$
i\eps \a_\eps u^\eps
+\eps^2/2 \, \D u^\eps
+n^2(x) u^\eps
=
1/\eps^{d/2}
S\(x/\eps\)
=:
S_\eps(x),
$
where for any function $f(x)$ we use the short-hand notation
$$
f_\eps(x)=
\frac{1}{\eps^{d/2}}
f\(\frac{x}{\eps}\)
.
$$
Using now the function $u^\eps$ instead of $w^\eps$, 
we observe the equality
\bea
\label{semeq}
&&
\<w^\eps,\phi\>
=
\<u^\eps,\phi_\eps\>
.
\eea
This transforms the original problem into the question of computing
the semiclassical limit $\eps\rgt 0$ in the equation satisfied by
$u^\eps$. One sees in (\ref{semeq}) that this limit
needs to be computed {\em at the semiclassical scale} (i.e. when tested upon
a smooth, concentrated function $\phi_\eps$).

In order to do so, we compute $u^\eps$ in terms of the 
semiclassical resolvent 
$\(
i\eps \a_\eps
+(\eps^2/2)\D
+n^2(x)
\)^{-1}
$.
It is the integral over the whole time interval
$[0,+\infty[$ of the propagator of the Schr\"odinger
operator associated with $\eps^2\D/2
+n^2(x)$. In other words we write
\bea
\non
u^\eps
&=&
\(
i\eps \a_\eps
+\frac{\eps^2}{2}\D
+n^2(x)
\)^{-1}
S_\eps
\\
&=&
i
\int_0^{+\infty}
\exp\(
i t \(
i\eps \a_\eps
+\frac{\eps^2}{2}\D
+n^2(x)
\)
\)
S_\eps
\,
dt
.
\eea
Now, defining the semi-classical propagator
\bea
U_\eps(t)
:=
\exp\(
i \frac{t}{\eps} \(
\frac{\eps^2}{2}\D
+n^2(x)
\)
\)
=
\exp\(
- i \frac{t}{\eps} H_\eps \)
,
\eea
associated with the semi-classical Schr\"odinger operator
\begin{align}
&&
H_\eps:=
-\frac{\eps^2}{2}\D
-n^2(x)
,
\end{align}
we arrive at the final formula
\bea
\label{integ}
\<w^\eps,\phi\>
=
\<u^\eps,\phi_\eps\>
=
\frac{i}{\eps}
\int_0^{+\infty}
\e^{-\a_\eps t} \;
\<
U_\eps(t)
S_\eps
,
\phi_\eps
\>
\, dt.
\eea
Our strategy is to pass to the limit in this very integral.

\bigskip

\noi
{\bf Second step: passing to the limit in the time integral (\ref{integ})}

In order to pass to the limit $\eps\rgt 0$ in (\ref{integ}),
we need to analyze the contributions
of various time scales in the corresponding time integral.
More precisely,
we choose for the whole subsequent analysis two (large) cutoff parameters
in time, denoted by $\T$ and $\tau$, and we analyze the contributions
to the time integral (\ref{integ}) that are due to the three regions
$$
0 \leq t \leq \T \eps
, \;
\T \eps \leq t \leq \tau
, \; \text{ and } \;
t \geq \tau . 
$$
We also choose a (small) exponent $\k>0$, and we occasionally
treat separately the contributions of very large times
$$
t \geq \eps^{-\k} .
$$
Associated with these truncations, we take once and for all
a smooth cutoff function $\chi$ defined on $\R$, such that
\bea
\label{chi}
&&
\non
\chi(z) \equiv 1 \; \text{ when } \; |z|\leq 1/2 ,
\; \;
\chi(z) \equiv 0  \; \text{ when } \; |z| \geq 1 ,
\\
&&
\qquad
\qquad
\qquad
\qquad
\chi(z)\geq 0 \; \text{ for any $z$}.
\eea
To be complete, there remains to finally choose a (small) cutoff parameter
in energy $\de>0$. Accordingly we distinguish in the $L^2$ scalar product
$\<U_\eps(t)S_\eps,\phi_\eps\>$ between energies close to (or far from) the
zero energy, which is critical for our problem. In other words, we set
the self-adjoint operator
$$
\chi_\de\(H_\eps\):=
\chi\(\frac{H_\eps}{\de}\) .
$$
This object is perfectly well defined using
standard functional calculus for self-adjoint operators.
We decompose
$$
\<U_\eps(t)S_\eps,\phi_\eps\>
=
\Big<U_\eps(t)\,\chi_\de(H_\eps)S_\eps,\phi_\eps\Big>
+
\Big<U_\eps(t)\, \( 1-\chi_\de\)(H_\eps)\,S_\eps,\phi_\eps\Big>
.
$$
Following the above described decomposition
of times and energies, we study each of the subsequent terms:

\medskip
\noi
$\bullet$ {\bf The contribution of small times} is
\beas
&&
\frac{1}{\eps}
\int_0^{2 \T \eps}
\chi\(\frac{t}{\T \eps}\)
\;
\e^{-\a_\eps t} \;
\<
U_\eps(t)
S_\eps
,
\phi_\eps\>
\; dt .
\eeas
We prove in section \ref{st} that this term actually gives the dominant
contribution in (\ref{integ}), provided the cutoff parameter
$\T$ is taken large enough. This (easy) analysis essentially boils down
to manipulations on the time dependent Schr\"odinger operator
$i \d_t + \D_x/2 + n^2(\eps x)$, for {\em finite times} $t$
of the order $t\sim\T$
at most.

\medskip
\noi
$\bullet$ {\bf The contribution of moderate and large
times, away from the zero energy},
is
\beas
&&
\frac{1}{\eps}
\int_{\T \eps}^{+\infty}
\(1-\chi\)\(\frac{t}{\T \eps}\)
\;
\e^{-\a_\eps t} \;
\Big<
U_\eps(t)
\(1-\chi_\de\)\(H_\eps\) S_\eps
,
\phi_\eps\Big>
\; dt .
\eeas
We prove in section \ref{2z}
below that this term has a vanishing contribution, provided
$\T$ is large enough. This easy result relies on a non-stationnary
phase argument in time, recalling
that $U_\eps(t)=\exp(-it H_\eps/\eps)$ and the energy $H_\eps$ is larger
than $\de>0$.

\medskip
\noi
$\bullet$ {\bf The contribution of very large times, close to the
zero energy} is
\beas
&&
\frac{1}{\eps}
\int_{\eps^{-\k}}^{+\infty}
\e^{-\a_\eps t} \;
\Big<
U_\eps(t)
\chi_\de\(H_\eps\) S_\eps
,
\phi_\eps
\Big>
\; dt .
\eeas
We prove in section \ref{kappa} that this term has a vanishing
contribution as $\eps\rgt 0$. To do so, we use
results proved by X.P. Wang \cite{Wa}: these essentially assert
that the operator $\<x\>^{-s} \; U_\eps(t) \chi_\de(H_\eps) \; \<x\>^{-s}$
has the natural size $\<t\>^{-s}$ as time goes to infinity,
provided the critical zero energy
is non-trapping.
Roughly, the semiclassical operator
$ U_\eps(t) \chi_\de(H_\eps)$ sends rays initially close to the origin,
at a distance of the order $t$ from the origin, when the energy is non
trapping. Hence the above scalar product involves
both a function $U_\eps(t)\chi_\de\(H_\eps\) S_\eps$ that is localized
at a distance $t$ from the origin, and a function $\phi_\eps$ that
is localized at the origin. This makes the corresponding contribution vanish.

\medskip

The most difficult terms are the last two that we describe now.

\medskip
\noi
$\bullet$ {\bf The contribution of large times, close to the
zero energy} is
\beas
&&
\frac{1}{\eps}
\int_{\tau}^{\eps^{-\k}}
\e^{-\a_\eps t} \;
\Big<
U_\eps(t)
\chi_\de\(H_\eps\) S_\eps,
\phi_\eps
\Big>
\; dt .
\eeas
The treatment of this term is performed in section \ref{brbr}.
It is similar in spirit to (though much harder than) the
analysis performed in the previous term: using only
information on the localization properties of
$U_\eps(t) \chi_\de\(H_\eps\) S_\eps$ and $\phi_\eps$,
we prove
that this term
has a vanishing contribution, provided $\tau$ is large enough.
To do so, we use ideas of Bouzouina and Robert \cite{BR}, to establish
a version of the Egorov theorem that holds true for {\em polynomially
large times} in
$\eps$.
We deduce that for any time
$\tau \leq t \leq \eps^{-\k}$, the term
$U_\eps(t)
\chi_\de\(H_\eps\) S_\eps$ is localized
close to the value at time $t$
of a trajectory shot from the origin.
The non-trapping assumption then says that for $\tau$ large enough,
$U_\eps(t)
\chi_\de\(H_\eps\) S_\eps$ is localized away from the origin.
This
makes the scalar product $\<U_\eps(t)
\chi_\de\(H_\eps\) S_\eps,\phi_\eps\>$
vanish asymptotically.

\medskip
\noi
$\bullet$ {\bf The contribution of moderate times close to the
zero energy} is
\beas
&&
\frac{1}{\eps}
\int_{\T \eps}^{\tau}
\(1-\chi\)\(\frac{t}{\T \eps}\)
\;
\e^{-\a_\eps t} \;
\Big<
U_\eps(t)
\chi_\de\(H_\eps\) S_\eps
,
\phi_\eps
\Big>
\; dt .
\eeas
This is the most difficult term: contrary to all preceding terms,
it cannot be analyzed using only geometric informations on the
microlocal support of the relevant functions.
Indeed, keeping in mind that the function
$U_\eps(t)
\chi_\de\(H_\eps\) S_\eps$ is localized on a trajectory initially shot from
the origin, whereas $\phi_\eps$ stays at the origin,
it is clear that
for times $\T \eps \leq t \leq \tau$,
the support of $U_\eps(t)
\chi_\de\(H_\eps\) S_\eps$ and $\phi_\eps$ {\em may intersect},
due to trajectories
passing {\em several times} at the origin. This might create a dangerous
accumulation of energy at this point.
For that reason, we need a precise
evaluation of the semi-classical propagator $U_\eps(t)$, for times
up to the order $t\sim\tau$. This is done using
the elegant wave-packet approach of M. Combescure and D. Robert
\cite{CRo} (see also
\cite{Ro}, and the nice lecture \cite{Ro2}): projecting $S_\eps$
over the standard gaussian wave packets, we can compute
$U_\eps(t) S_\eps$
in a quite explicit fashion, with the help of
classical quantities like, typically, the linearized flow of the Hamiltonian
$\xi^2/2-n^2(x)$. This gives us an integral representation
with a complex valued phase function. Then,
one needs to insert a last (small) cutoff parameter in time,
denoted $\th>0$.
For small times,
using the above mentioned representation formula, we first prove 
that the term
\beas
&&
\frac{1}{\eps}
\int_{\T \eps}^{\th}
\(1-\chi\)\(\frac{t}{\T \eps}\)
\;
\e^{-\a_\eps t} \;
\Big<
U_\eps(t)
\chi_\de\(H_\eps\) S_\eps
,
\phi_\eps
\Big>
\; dt ,
\eeas
vanishes asymptotically,
provided $\th$ is small, and $\T$ is large enough.
To do so, we use that for small enough $\th$,
the propagator $U_\eps(t)$ acting on $S_\eps$ resembles the free
Schr\"odinger operator $\exp\(i t [\D_x/2+n^2(0)]\)$. 
In terms of trajectories, on this time scale,
we use that $U_\eps(t) S_\eps$ is localized around a ray
that leaves the origin {\em at speed $n(0)$}.
Then, for later times, we prove
that the remaining contribution
\beas
&&
\frac{1}{\eps}
\int_{\th}^{\tau}
\e^{-\a_\eps t} \;
\Big<
U_\eps(t)
\chi_\de\(H_\eps\) S_\eps
,
\phi_\eps
\Big>
\; dt ,
\eeas
is small. This uses stationary
phase formulae in the spirit of \cite{CRR}, and this is where
the transversality assumption {\bf (H)} page \pageref{HH} enters:
trajectories passing
several times at the origin do not accumulate to much
energy at this point.

We end up this sketch of proof with a figure illustrating the
typical trajectory (and the associated cutoffs in time) that 
our analysis has to deal with.
\begin{center}
\scalebox{0.42}{\includegraphics{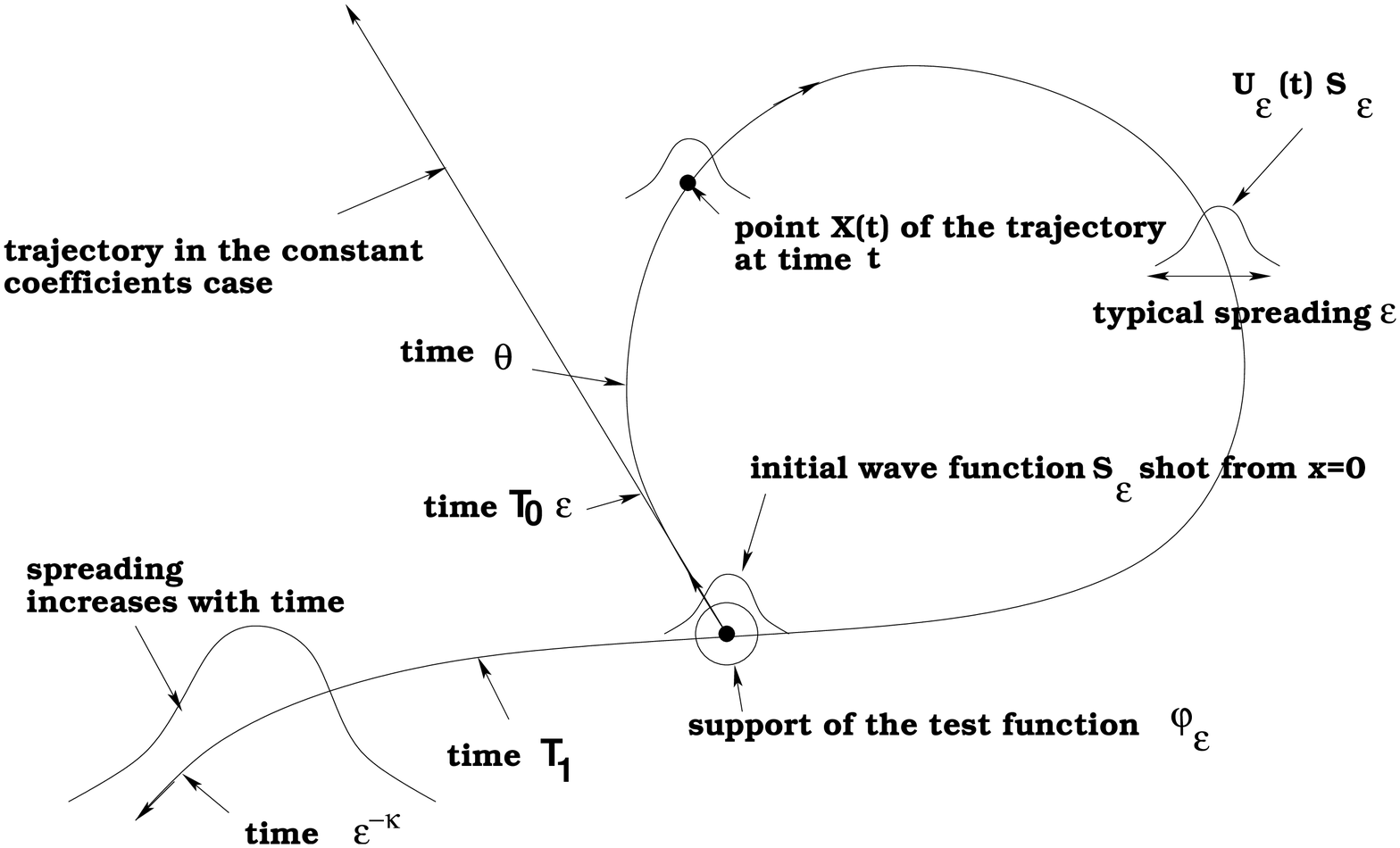}}
\end{center}



\subsection{Notations used in the proof}


Throughout this article, we will make use of the following notations.

\medskip
\noi
$\bullet$
{\bf Semi-classical quantities}

\noi
The semi-classical Hamiltonian $H_\eps$ and its associated
propagator $U_\eps(t)$ have already been defined.
We also need to use the Weyl quantization.
For a symbol $a(x,\xi)$ defined on $\R^{2d}$,
its Weyl quantization is
\beas
\({\rm Op}^w_\eps(a) f\)(x)
:=
\frac{1}{(2\pi \eps)^d}
\int_{\R^{2d}}
\e^{i\frac{(x-y) \cdot \xi}{\eps}}
\;
a\(\frac{x+y}{2},\xi\)
\;
f(y) \; dy \; d\xi .
\eeas
Throughout the paper, we use the standard semi-classical symbolic calculus,
and refer, e.g., to \cite{DS} or \cite{Ma}. In particular,
for a weight $m(x,\xi)$, we use symbols $a(x,\xi)$ in the class
$S(m)$, i.e. symbols such that for any multi-index $\a$, there
exists a constant $C_\a$ so that 
$$
|\d^\a a(x,\xi)| \leq C_\a m(x,\xi) , \; \forall (x,\xi) \in \R^{2d} .
$$
The notation $a\sim\sum\eps^k a_k$ means that for any $N$ and any $\a$,
there exists a constant $C_{N,\a}$ such that
$$
\Bigg|
\d^\a \(a(x,\xi)-\sum_{k=0}^N \eps^k a_k(x,\xi)\)
\Bigg| \leq C_{N,\a} \; \eps^{N+1}  m(x,\xi) ,
\; \forall (x,\xi) \in \R^{2d} .
$$

\medskip

\noi
$\bullet$
{\bf Classical quantitities}

\noi
Associated with the Hamiltonian
$H(x,\xi)=\xi^2/2-n^2(x)$,
we denote the Hamiltonian flow
\beas
&&
\Phi(t,x,\xi)=(X(t,x,\xi)\, , \, \Xi(t,x,\xi)) ,
\eeas
defined as the solution of the Hamilton equations
\begin{align}
\label{hamf}
\non
&\disp\frac{\d}{\d t} X(t,x,\xi)
=
\Xi(t,x,\xi) 
,
&
X(0,x,\xi)=x ,
\qquad\qquad
\\
&\disp\frac{\d}{\d t} \Xi(t,x,\xi)
=
\(\nabla_x n^2\)\(X(t,x,\xi)\) 
,
&
\Xi(0,x,\xi)=\xi .
\qquad\qquad
\end{align}
These may be written shortly
\bea
\label{hamfs}
\frac{\d}{\d t} \Phi(t,x,\xi)
=
J \; \frac{D H}{D(x,\xi)}\(\Phi(t,x,\xi)\)
,
\eea
where $J$ is the standard symplectic matrix
\bea
\label{j}
J=
\left(
\begin{array}{cc}
0&{\rm Id}
\\
-{\rm Id}&0
\end{array}
\right) .
\eea
The linearized flow of $\Phi$ is denoted by
\bea
\label{linf}
F(t,x,\xi):=\frac{D \Phi(t,x,\xi)}{D(x,\xi)} .
\eea
It may be decomposed into
\bea
\label{decf}
F(t,x,\xi)
=
\left(
\begin{array}{cc}
A(t,x,\xi)&B(t,x,\xi)
\\
C(t,x,\xi)&D(t,x,\xi)
\end{array}
\right)
,
\eea
where the matrices $A(t)$, $B(t)$, $C(t)$, and $D(t)$ are, by definition
\beas
&&
A(t,x,\xi)=\frac{D X(t,x,\xi)}{D x} 
, \quad
B(t,x,\xi)
=\frac{D X(t,x,\xi)}{D \xi}
,
\\
&&
C(t,x,\xi)
=
\frac{D \Xi(t,x,\xi)}{D x}
, \quad
D(t,x,\xi)
=
\frac{D \Xi(t,x,\xi)}{D \xi} .
\eeas
Upon linearizing (\ref{hamf}),
the matrices $A(t)$, $B(t)$, $C(t)$, and $D(t)$ clearly satisfy
the differential system
\begin{align}
\label{linf1}
\non
&
\frac{\d}{\d t} A(t,x,\xi)=C(t,x,\xi) ,
&
A(0,x,\xi)={\rm Id} ,
\\
&
\frac{\d}{\d t} C(t,x,\xi)=\frac{D^2 n^2}{D x^2}
\(X(t,x,\xi)\) \; A(t,x,\xi) ,
&
C(0,x,\xi)=0 ,
\end{align}
together with
\begin{align}
\label{linf2}
\non
&
\frac{\d}{\d t} B(t,x,\xi)=D(t,x,\xi) ,
&
B(0,x,\xi)=0 ,
\\
&
\frac{\d}{\d t} D(t,x,\xi)=\frac{D^2 n^2}{D x^2}\(X(t,x,\xi)\) \; B(t,x,\xi)
,
&
D(0,x,\xi)={\rm Id} .
\end{align}
In short, one may write as well
\bea
\label{hamlinfs}
\frac{\d}{\d t} F(t,x,\xi)
=
J \; \frac{D^2 H}{D(x,\xi)^2}\(\Phi(t,x,\xi)\) \; F(t,x,\xi) .
\eea
A last remark is in order. Indeed, it is a standard fact to observe that
the matrix $F(t,x,\xi)$ is a symplectic matrix, in that
\bea
F(t,x,\xi)^\tr J F(t,x,\xi)=J ,
\eea
for any
$(t,x,\xi)$. Here, the exponent $\tr$ denotes
transposition. Decomposing $F(t)$ as in (\ref{decf}),
this gives the relations
\bea
\label{relsym}
\non
&&
A(t)^\tr C(t)=C(t)^\tr A(t) , \;
B(t)^\tr D(t)=D(t)^\tr B(t) , \;
\\
&&
A(t)^\tr D(t)-C(t)^\tr B(t)={\rm Id} .
\eea
These can be put in the following useful form
\begin{align}
\label{uv}
\non
&
\(A(t)+iB(t)\)^\tr \(C(t)+iD(t)\)
=
\(C(t)+iD(t)\)^\tr \(A(t)+iB(t)\)
&
\\
&
\(C(t)+iD(t)\)^\tr \(A(t)-iB(t)\)
&
\\
&
\non
\qquad
\qquad
-
\(A(t)+iB(t)\)^\tr \(C(t)-iD(t)\)
=
2 i {\rm Id} .
&
\end{align}
These relations will be used in section \ref{mod}.


\section{Small time contribution: the case
$0 \leq t \leq \T \eps$}
\label{st}

\setcounter{equation}{0}


In this section, we prove the following
\begin{proposition}
\label{ouane}
We use the notations of section \ref{redu}.
The refraction index $n^2$ is assumed bounded and continuous.
The data $S$ and $\phi$ are supposed to belong to ${\cal S}(\R^d)$.
Then, the following holds:
\vskip0.1cm
\noi
(i)
for any fixed value
of $\T$, we have the asymptotics
\bea
\non
&&
\frac{i}{\eps}
\int_0^{2 \T \eps}
\chi\(\frac{t}{\T \eps}\)
\;
\e^{-\a_\eps t} \;
\<
U_\eps(t)
S_\eps
,
\phi_\eps\>
\; dt
\\
&&
\qquad
\mathop{\longrightarrow}\limits_{\eps\rgt 0}
i
\int_0^{2\T}
\chi\(\frac{t}{\T}\)
\;
\<
\exp\(i t (\D_x/2+n^2(0))\)
\;
S
,
\phi\>
\; dt
.
\label{assmal}
\eea
\vskip0.1cm
\noi
(ii)
Besides, there exists a 
universal constant $C_d$ depending only on the dimension,
such that the right-hand-side of (\ref{assmal}) satisfies
\bea
\label{t0inf}
\non
&&
\Bigg|
i
\int_0^{2\T}
\chi\(\frac{t}{\T}\)
\;
\<
\exp\(i t (\D_x/2+n^2(0))\)
\;
S
,
\phi\>
\; dt
-
\<w^{\rm out},\phi\>
\Bigg|
\\
&&
\qquad\qquad\qquad
\leq
C_d \; \T^{-d/2+1}
\mathop{\longrightarrow}\limits_{\T\rgt \infty}
0
.
\eea
\end{proposition}

\noi
{\bf Proof of proposition \ref{ouane}}

\noi
{\it Part (i)}

\noi
In order to recover the limiting value announced in (\ref{assmal}),
we first perform the inverse scaling that leads from $w^\eps$ to
$u^\eps$ (see (\ref{wtou})). We rescale time
$t$ by a factor $\eps$ as well. This gives
\beas
&&
\frac{1}{\eps}
\int_0^{+\infty}
\chi\(\frac{t}{\T \eps}\)
\;
\e^{-\a_\eps t} \;
\<
U_\eps(t)
S_\eps
,
\phi_\eps
\>
\; dt
\\
&&
\qquad
=
\int_0^{+\infty}
\chi\(\frac{t}{\T}\)
\;
\e^{-\eps \a_\eps t} \;
\<
U_\eps(\eps \, t)
S_\eps
,
\phi_\eps
\>
\; dt
\\
&&
\qquad
=
\int_0^{+\infty}
\chi\(\frac{t}{\T}\)
\;
\e^{-\eps \a_\eps t} \;
\<
\exp\(i t \(\D/2 + n^2(\eps x)\)\)
S
,
\phi
\>
\; dt .
\eeas
We now let
\beas
&&
{\bf w}^\eps(t,x):=
\exp\(i t \(\D/2 + n^2(\eps x)\)\)
S(x)
.
\eeas
The function ${\bf w}^\eps(t,x)$ is bounded
in $L^\infty\(\R;L^2\(\R^d\) \)$, and it satisfies in the
distribution
sense
\beas
&&
i \d_t {\bf w}^\eps(t,x)
=
-\frac12\D_x {\bf w}^\eps(t,x)
- n^2(\eps x) {\bf w}^\eps
,
\quad
{\bf w}^\eps(0,x)
=
S(x) .
\eeas
These informations are enough to deduce that there exists
a function 
$
{\bf w}(t,x) \in L^\infty\(\R;L^2\(\R^d\)\)
$
such that
a subsequence of  ${\bf w}^\eps(t,x)$ goes,
as $\eps \rgt 0$,  to ${\bf w}(t,x)$ in 
$L^\infty\(\R;L^2\(\R^d\)\)$ - weak$*$. On the more, the
limit ${\bf w}(t,x)$ obviously
satisfies in the distribution sense
\beas
&&
i \d_t {\bf w}(t,x)
=
-\frac12\D_x {\bf w}(t,x)
- n^2(0) {\bf w}
, \quad
{\bf w}(0,x)
=
S(x) .
\eeas
In other words
$$
{\bf w}(t)
=\exp\(i t \(\D/2 + n^2(0)\)\)
S(x) .
$$
Hence, by uniqueness of the limit, the whole sequence
${\bf w}^\eps(t,x)$ goes to ${\bf w}(t,x)$
in $L^\infty\(\R;L^2\(\R^d\)\)$-weak$*$. This proves (\ref{assmal})
and part (i) of the proposition.

\medskip

\noi
{\it Part (ii)}

\noi
This part is easy and relies on the standard dispersive properties of the
free Schr\"odinger equation.
Indeed, we have
\begin{align*}
&
\Big|
\<\exp\(i t \(\D_x/2+n^2(0)\) \) S, \phi \>
\Big|
&
\\
&
\qquad\qquad
\leq
\Big\|\exp\(i t \(\D_x/2+n^2(0)\) \) S \,\Big\|_{L^\infty} \; \|\phi\|_{L^1}
&
\\
&
\qquad\qquad
\leq
C_d \; t^{-d/2} \;  \|S\|_{L^1} \; \|\phi\|_{L^1}
,
&
\end{align*}
(recall that $S$ and $\phi$ are assumed smooth enough to have
finite $L^1$ norm), for some constant $C_d>0$ that only depends
upon the dimension $d$. This, together with the integrability
of the function 
$ t^{-d/2}$ at infinity when $d\geq 3$, ends the proof
of (\ref{t0inf}).
\qed


\section{Contribution of mo\-de\-ra\-te and large times,
a\-way from the zero energy}
\label{2z}

\setcounter{equation}{0}


In this section we prove the (easy)
\begin{proposition}
\label{2}
We use the notations of section \ref{redu}.
The index $n^2$ is assumed to have the symbolic behaviour
(\ref{decn}). The data $S$ and $\phi$ are supposed to belong to $L^2(\R^d)$.
Then, there exists a constant $C_\de>0$,
which depends
on the cutoff parameter
$\de$, such that for any $\eps \leq 1$, and $\T\geq 1$, we have
\bea
\non
&&
\label{majaw}
\Bigg|
\frac{1}{\eps}
\int_{\T \eps}^{+\infty}
\(1-\chi\)\(\frac{t}{\T \eps}\)
\;
\e^{-\a_\eps t} \;
\Big<
\(1-\chi_\de\(H_\eps\)\)U_\eps(t)
S_\eps
,
\phi_\eps
\Big>
\; dt \;
\Bigg|
\\
&&
\qquad\qquad\qquad
\leq C_\de \(\frac{1}{\T}+\a_\eps^2 \)
.
\eea
\end{proposition}

\medskip

\noi
{\bf Proof of proposition \ref{2}}

\noi
The proof relies on a simple non-stationary phase argument. Indeed,
this term has the value
$$
\frac{1}{\eps}
\int_0^{+\infty}
\(1-\chi\)\(\frac{t}{\T \eps}\)
\;
\e^{-\a_\eps t} \;
\Big<
\(1-\chi_\de\(H_\eps\)\)\exp\(-i \frac{t}{\eps} H_\eps\)
S_\eps
,
\phi_\eps
\Big>
\; dt.
$$
Hence, making the natural integrations by parts in time,
we recover the value
\begin{align*}
&
\eps^2
\int_0^{+\infty}
\frac{\d^3}{\d t^3}
\(
\(1-\chi\)\(\frac{t}{\T \eps}\)
\e^{-\a_\eps t} \,
\)
\\
&
\qquad\qquad\qquad\qquad
\Bigg<
\frac{\(1-\chi_\de\(H_\eps\)\)}{(-i H_\eps)^3}
\exp\(-i \frac{t}{\eps} H_\eps\)
S_\eps
,
\phi_\eps
\Bigg>
\, dt .
\end{align*}
A direct inspection shows that this is bounded by
\beas
&&
C \; \eps^2 \; \de^{-3} \; \|S\|_{L^2} \; \|\phi\|_{L^2}
\;
\int_0^{+\infty}
\Bigg|
\frac{\d^3}{\d t^3}
\(
\(1-\chi\)\(\frac{t}{\T \eps}\)
\e^{-\a_\eps t} \;
\)
\Bigg|
\, dt
\\
&&
\qquad
\leq
C \; \eps^2 \; \de^{-3} \;
\;
\|\chi\|_{W^{3,\infty}}
\;
\(
\frac{1}{\T^2 \eps^2}+
\frac{1}{\T \eps}+
\a_\eps^2+\a_\eps^2\) .
\eeas
\qed


\section{Contribution of large times, close to the zero energy: the
case $t \geq \eps^{-\k}$}
\label{kappa}

\setcounter{equation}{0}


In this section we prove the following
\begin{proposition}
\label{wawang}
We use the notations of section \ref{redu}.
The index $n^2$ is assumed to have the symbolic behaviour (\ref{decn}).
The Hamiltonian flow associated with
$\xi^2/2-n^2(x)$ is assumed non-trapping at the zero energy level.
Finally, the data $S$ and $\phi$ are supposed to belong
to ${\cal S}(\R^d)$.
Then, for any $\de>0$ small enough, and for any
$\k>0$,
there exists a constant $C_{\k,\de}$ depending on $\k$ and $\de$, so that
\bea
\label{lar}
&&
\Bigg|
\frac{1}{\eps}
\int_{\eps^{-\k}}^{+\infty}
\e^{-\a_\eps t} \;
\Big<
U_\eps(t)
\chi_\de\(H_\eps\) S_\eps
,
\phi_\eps
\Big>
\; dt \;
\Bigg|
\leq
C_{\k,\de} \; \eps
.
\eea
\end{proposition}

\medskip

The proof relies on the dispersive properties
of the semi-classical propagator
$U_\eps(t)$, inherited from the ones of the classical
flow $\Phi(t)$.
More quantitatively, we use in this section
a Theorem by X.P. Wang \cite{Wa}, that we now state.
Our index of refraction $n^2(x)$ is such that
$n^2(x)$ lies in $C^\infty(\R^d)$, and it has the symbolic behaviour
\beas
&&
n^2(x)=n^2_\infty - V(x)
, \; \text{ with }
|\d^\a V(x)| \leq \<x\>^{-\rho-|\a|}
\eeas
(the case $0<\rho\leq 1$ is the long-range case, and
the case $\rho > 1$ is the short-range case, in the terminology
of quantum scattering). On the more, the trajectories of the classical
flow at the zero energy (i.e. on the set
$\{(x,\xi) \in \R^{2d} \; \text{ s.t. } \; \xi^2/2-n^2(x)=0 \}$)
are assumed non-trapped. It is known \cite{DG} that this non-trapping
behaviour
is actually an open property, in that
\bea
\label{opennontrap}
\non
&&
\text{
there exists a $\de_0>0$ such that
for any energy $E$
}
\\
&&
\text{
satisfying $|E|\leq \de_0$, the trajectories
of the classical flow}
\\
\non
&&
\text{
at the energy $E$ are non-trapping as well.
}
\eea
Under these circumstances,
it has been proved in \cite{Wa} that for any
real $s>0$, and for any $\eta>0$, the following weighted estimate
holds true,
\bea
\label{wang}
\forall t \in \R ,
\quad
\|\<x\>^{-s} U_\eps(t) \chi_\de(H_\eps) f\|_{L^2}
\leq
\frac{C_{\de,\eta,s}}{\<t\>^{s-\eta}}
\; \; 
\|\<x\>^s f(x)\|_{L^2} ,
\eea
provided the cutoff
in energy $\de$ satisfies $\de\leq\de_0$, i.e. provided we are only
looking at trajectories having a non-trapping energy.
This inequality holds
for any test function $f$, and for some constant $C_{\de,\eta,s}$
depending only on $\de$, $\eta$ and $s$.
In the short-range case ($\rho>1$), one may even take $\eta=0$ in the above
estimate.
Note that
\cite{Wa}
actually proves more: in some sense,
the non-trapping behaviour of the classical flow is {\em equivalent}
to the time decay (\ref{wang}). We refer to the original
article for details. We are now ready to give the

\medskip

\noi
{\bf Proof of proposition \ref{wawang}}

\noi
Taking 
$\de\leq\de_0$, we estimate, using (\ref{wang}),
\beas
&&
\frac{1}{\eps}
\Bigg|
\int_{\eps^{-\k}}^{+\infty}
\;
\e^{-\a_\eps t} \;
\<
\chi_\de\(H_\eps\)
U_\eps(t)
S_\eps
,
\phi_\eps
\>
\; dt
\Bigg|
\\
&&
\qquad
\leq
\frac{1}{\eps}
\int_{\eps^{-\k}}^{+\infty}
\;
\|\<x\>^{-s} U_\eps(t) \chi_\de(H_\eps) S_\eps\|_{L^2}
\;
\|\<x\>^s \phi_\eps\|_{L^2}
\; dt
\\
&&
\qquad
\leq
\frac{1}{\eps}
\;
\|\<x\>^s  S_\eps(x)\|_{L^2}
\;
\;
\|\<x\>^s \phi_\eps\|_{L^2}
\int_{\eps^{-\k}}^{+\infty}
\;
\frac{C_{\de,\eta,s}}{\<t\>^{s-\eta}}
\; dt
\\
&&
\qquad
\leq
C_{\de,\eta,s} \; 
\eps^{\k(s-\eta-1)-1}
\;
\|\<x\>^s  S_\eps(x)\|_{L^2}
\;
\|\<x\>^s \phi_\eps\|_{L^2}
\eeas
Hence, taking $s$ large enough, and $\eta$ small enough, e.g.
$s=2+2/\k$, $\eta=1$,
we obtain an upper bound of the size
$$
C_{\k,\de} \; \eps \; \|\<x\>^s  S(x)\|_{L^2}
\;
\|\<x\>^s \phi\|_{L^2} .
$$
Here we used the easy fact that
$\|\<x\>^s  f_\eps(x)\|_{L^2}\leq
\|\<x\>^s  f(x)\|_{L^2}$, when $\eps \leq 1$, together
with $\|\<x\>^s  S(x)\|_{L^2}<\infty$, and similarly for $\phi$.
\qed


\section{Contribution of large times, close to the zero energy:
the case $\tau \leq t \leq \eps^{-\k}$}
\label{brbr}

\setcounter{equation}{0}


To complete the analysis of the contribution
of ``large times'' and ``small energies'' in (\ref{integ}) that we begun
in section \ref{kappa}, there
remains to estimate the term
\bea
\label{rem3}
&&
\frac{1}{\eps}
\int_{\tau}^{\eps^{-\k}}
\(1-\chi\)\(\frac{t}{\tau}\)
\;
\e^{-\a_\eps t} \;
\<
\chi_\de\(H_\eps\)
U_\eps(t)
S_\eps
,
\phi_\eps
\>
\; dt.
\eea
In this section, we prove the
\begin{proposition}
\label{brob}
We use the notations of section \ref{redu}.
The index $n^2$ is assumed to have the symbolic behaviour
(\ref{decn}) with $n_\infty^2>0$\footnote{
The assumption $n_\infty^2$ is crucial, see Lemma \ref{bourob} below.
It ensures that the wave $U_\eps(t) S_\eps$ propagates
with a uniformly non-zero speed, at infinity in time $t$.
}.
The Hamiltonian flow associated with
$\xi^2/2-n^2(x)$ is assumed non-trapping at the zero energy.
Finally, the data $S$ and $\phi$ are supposed to belong to
${\cal S}(\R^d)$. Then,
for $\de>0$ small enough,
there exists a $\tau(\de)$ depending
on $\de$ such that
for any $\tau\geq\tau(\de)$, we have for $\k$ small enough,
\bea
\label{rem3est}
\non
&&
\Bigg|
\frac{1}{\eps}
\int_{\tau}^{\eps^{-\k}}
\(1-\chi\)\(\frac{t}{\tau}\)
\;
\e^{-\a_\eps t} \;
\<
\chi_\de\(H_\eps\)
U_\eps(t)
S_\eps
,
\phi_\eps
\>
\; dt \;
\Bigg|
\\
&&
\qquad\qquad\qquad
\leq
C_{\k,\de} \; \eps , \; \text{ as } \; \eps \rgt 0 ,
\eea
for some constant $C_{\k,\de}$ that depends upon $\k$ and $\de$.
\end{proposition}

\medskip

The idea of proof is the following:
the functions $S_\eps$ and $\phi_\eps$ are
microlocally supported close to points $(x_0,\xi_0) \in \R^{2d}$
such that $x_0=0$ (due to the concentration
of both functions close to the origin as $\eps\rgt 0$).
On the more, using the Egorov Theorem, one may think
of the time evolved function
$U_\eps(t)
S_\eps$ as being microlocally supported
close to points $(X(t;x_0,\xi_0),\Xi(t;x_0,\xi_0))$ that are
trajectories of the classical flow, with initial data
$(x_0,\xi_0)$ such that $x_0=0$.
Using the non-trapping assumption on the classical flow,
we see that for large times $t\geq \tau$ with $\tau$ large enough,
the trajectory $X(t;x_0,\xi_0)$ with $x_0=0$ is far away
from the origin. Hence the microlocal support
of $U_\eps(t)
S_\eps$  and $\phi_\eps$ do not intersect, and the
factor (\ref{rem3}) should be arbitrary small in $\eps$ as $\eps \rgt 0$.

The difficulty in making this last statement rigorous
lies in the fact that we need to use the Egorov Theorem up to 
(polynomially) large
times of the order $t\sim \eps^{-\k}$. This difficulty
is solved in Lemma \ref{bourob} below. Indeed, upon
adapting a recent result of Bouzouina and Robert \cite{BR}
we give remainder estimates in the Egorov Theorem that hold
up to polynomially
large times (logarithmic times are obtained in the context of \cite{BR}).
This is enough to conclude.


\subsection{Proof of proposition \ref{brob}}


The proof is given in several steps.

\medskip

\noi
{\it {\bf First step}: Preliminary reduction}

In this step we quantify the fact that the functions involved in the scalar
product in (\ref{rem3est}) are microlocalized close to the zero energy
$\xi^2/2=n^2(x)$ (in frequency) and close to the origin $x=0$ (in space).
To do so, we simply write, using the fact that $S$ and $\phi$
belong to ${\cal S}(\R^d)$,
\beas
&&
\phi_\eps(x) = \chi_\de(|x|) \phi_\eps(x)
+O_\de(\eps^\infty) \quad \text{ in } L^2(\R^d),
\eeas
and similarly for $S_\eps$.
This means that for any integer $N$, there exists a $C_{N,\de}>0$ that depends
on $N$ and $\de$, such that
$\|\phi_\eps(x) - \chi_\de(|x|) \phi_\eps(x)\|_{L^2(\R^d)} \leq C_N \eps^N$.
As a consequence, we may rewrite the contribution (\ref{rem3})
we are interested in as
\beas
\frac{1}{\eps}
\int_{\tau}^{\eps^{-\k}}
\(1-\chi\)\(\frac{t}{\tau}\)
\;
\e^{-\a_\eps t} \;
\<
\chi_\de(|x|) \,
\chi_\de\(H_\eps\)
U_\eps(t)
\chi_\de(|x|) \,
S_\eps
,
\phi_\eps
\>
\; dt
\eeas
up to an $O_\de(\eps^\infty)$. There remains to bound the above term by
\bea
\label{32}
\non
&&
\leq
\|S_\eps\|_{L^2}
\,
\|\phi_\eps\|_{L^2}
\times
\frac{1}{\eps}
\int_{\tau}^{\eps^{-\k}}
\,
\Big\|
\chi_\de\(|x|\)
\chi_\de\(H_\eps\)
U_\eps(t)
\chi_\de\(|x|\)
\Big\|_{{\cal L}(L^2)}
\, dt 
\\
&&
\leq
\frac{C}{\eps}
\,
\int_{\tau}^{\eps^{-\k}}
\,
\Big\|
\chi_\de\(|x|\)
\chi_\de\(H_\eps\)
U_\eps(t)
\chi_\de\(|x|\)
\Big\|_{{\cal L}(L^2)}
\, dt 
\, ,
\eea
up to an $O_\de(\eps^\infty)$.
Our strategy is to now evaluate the operator norm under the integral sign.
This task is performed in the next two steps.

\bigskip

\noi
{\it {\bf Second step:}
symbolic calculus}

\noi
In view of (\ref{32}), our analysis boils down to computing,
for any $\tau \leq t \leq \eps^{-\k}$, the operator norm
\begin{align*}
&&
\Big\|
\chi_\de\(|x|\)
\chi_\de\(H_\eps\)
U_\eps(t)
\chi_\de\(|x|\)
\Big\|^2_{{\cal L}(L^2)}.
\end{align*}
Expanding the square,
this norm has the value
\begin{align}
\label{devas}
&&
\Big\|
\chi_\de\(|x|\)
U^*_\eps(t)
\chi_\de\(H_\eps\)
\chi_\de^2\(|x|\)
\chi_\de\(H_\eps\)
U_\eps(t)
\chi_\de\(|x|\)
\Big\|_{{\cal L}(L^2)}.
\end{align}
Now, and for later convenience, we rewrite the above localizations
in energy and space, as microlocalisations in position and frequency.

Using the functional calculus for pseudodifferential operators
of Helffer and Robert
\cite{HR} (see also the lecture notes \cite{DS} and \cite{Ma}),
there exists
a symbol $\Chid(x,\xi)$
such that
$$
\chi_\de\(H_\eps\)
=
{\rm Op}^w_\eps(\Chid)+O(\eps^\infty)
\quad
\text{ in }
\; {\cal L}(L^2) .
$$
The symbol $\Chid(x,\xi)$ is given by a formal expansion
\bea
\label{hchi}
&&
\Chid(x,\xi)
\sim
\sum_{k\geq 0}
\eps^k
\Chid^{(k)}(x,\xi),
\eea
where the expansion (\ref{hchi}) holds in the class of symbols
that are bounded together with all their derivatives.
Furthermore, the principal symbol of $\Chid$ is computed through the natural
equality
$$
\Chid^{(0)}(x,\xi)=\chi_\de\(\frac{\xi^2}{2}-n^2(x)\) .
$$
Finally, the explicit formulae in \cite{DS}
give at any order $k\geq 0$ the following information on the support
of the symbols $\Chid^{(k)}$,
$$
{\rm supp} \; \Chid^{(k)} \subset 
\{ |\xi^2/2-n^2(x) |\leq \de \} .
$$
Hence (\ref{devas}) becomes, using standard symbolic calculus,
\bea
\label{34}
\Big\|
\chi_\de\(|x|\)
\,
U^*_\eps(t)
\,
\[
{\rm Op}_\eps^w(\Chid(x,\xi) \,\,\sharp\,
\chi_\de^2\(|x|\) \,\sharp\,  \Chid(x,\xi) ) 
\]
\,
U_\eps(t)
\,
\chi_\de\(|x|\)
\Big\|_{{\cal L}(L^2)},
\eea
up to an $O_\de(\eps^\infty)$ (Here we used the uniform bound
$\|U_\eps(t)\|_{{\cal L}(L^2)}\leq 1$). Let us define for convenience
the following short-hand notation for the symbol in brackets in (\ref{34}):
$$
\bd(x,\xi)
:=
\Chid(x,\xi) \,\,\sharp\,
\chi_\de^2\(|x|\) \,\sharp\,  \Chid(x,\xi).
$$
The only information we need in the sequel is that $\bd$ admits an asymptotic
expansion $\bd=\sum_{k\geq 0} \eps^k b_\de^{(k)}$, where each
$b_\de^{(k)}$ has support
\beas
\non
{\rm supp} \; \bd^{(k)}
&\subset&
\{
|x|\leq \de \}
\cap
\{|\xi^2/2-n^2(x)| \leq \de \}
=:E(\de).
\eeas
This serves as a definition of the
(compact) set $E(\de)$ in phase space.
In the sequel, we summarize these informations
in the following
abuse of notation
\bea
\label{supb}
&&
{\rm supp} \; \bd
\subset
E(\de).
\eea
The remainder part of our analysis is devoted to estimating
$$
\Big\|
\chi_\de\(|x|\)
\,
U^*_\eps(t)
\,
{\rm Op}_\eps^w(
\bd(x,\xi))
\,
U_\eps(t)
\,
\chi_\de\(|x|\)
\Big\|_{{\cal L}(L^2)},
$$
and the hard part of the proof lies in establishing an
``Egorov theorem for large times'', to compute the
conjugation $U^*_\eps(t) {\rm Op}_\eps^w(
\bd(x,\xi)) U_\eps(t)$ in (\ref{devas}).

\bigskip

\noi
{\it {\bf Third step:}
an Egorov theorem valid for large times - End of the proof}

\noi
Now we claim the following
\begin{lemma}
\label{bourob}
\noi
We assume that the refraction index
has the symbolic behaviour (\ref{decn}) with $n_\infty^2>0$\footnote{
The assumption $n_\infty^2>0$ is crucial, see (\ref{cde})}.
We also assume that
the zero energy is non-trapping for the flow.
Take the cutoff parameter in energy  $\de$ small enough.
Then,

(i)
Let $\Phi(t,x,\xi)$ be the classical flow associated with the
Hamiltonian $\xi^2/2-n^2(x)$. 
Let $F(t,x,\xi)$ be the linearized flow. For any
multi-index $\a$, and for any (small) parameter $\eta>0$,
there exists a constant $C_{\de,|\a|,\eta}$ such that
for any initial datum $(x,\xi) \in E(\de)
=\{
|x|\leq \de \}
\cap
\{|\xi^2/2-n^2(x)| \leq \de \}
$, we have
\bea
\label{linflow}
\forall t \in \R , \quad
\Bigg|
\frac{ \d^\a F(t,x,\xi)}{\d (x,\xi)^\a}
\Bigg|
\leq
C_{\de,|\a|,\eta} \;
\<t\>^{(1+\eta)(1+|\a|)+2|\a|}.
\eea
In other words, the linearized flow has at most polynomial growth with time.

\medskip

\noi
(ii)
As a consequence,
for any time $t$, there exists a time-dependent symbol
\beas
\bbd(t,x,\xi)
\sim
\sum_{k\geq 0}
\eps^k
\bbd^{(k)}(t,x,\xi)
,
\eeas
such that the following holds: there exists a number $c_\de>0$ such that
for any $N>0$, there exists a constant
$C_{\de,N}$ such that
\begin{align}
\label{remestbr}
&
\hspace{-0.6cm}
\Bigg\|
U^*_\eps(t)
{\rm Op}^w_\eps\( \bd \)
U_\eps(t)
-
{\rm Op}_\eps^w\(
\sum_{k=0}^N \eps^k \bbd^{(k)}
\)
\Bigg\|_{{\cal L}(L^2)}
\leq
C_{\de,N} \,
\eps^{N+1} \,
\<t\>^{c_\de N^2}.
\end{align}
Again, the error grows polynomially with time,
and we have some control on the dependence
of the estimates with the truncation parameter $N$.

\medskip

\noi
(iii)
Moreover, we have the natural formulae
\beas
&&
\bbd^{(0)}(t,x,\xi)
=
\bd\(\Phi(t,x,\xi)\)
,
\eeas
and, for any $k\geq 0$ we have the information on the support
\beas
{\rm supp} \;  \bbd^{(k)}(t,x,\xi)
\subset
\{(x,\xi)\in \R^{2d} \; \text{ s.t. } \;
\Phi(t,x,\xi) \in E(\de) \} .
\eeas
\end{lemma}
\medskip

\noi
We postpone the proof of Lemma \ref{bourob} to paragraph
\ref{6.2} below.
We first draw its consequences in our perspective.

\medskip

Leaving $N$ as a free parameter for the moment,
we obtain
\begin{align*}
&
\Big\|
\chi_\de\(|x|\)
\,
U^*_\eps(t)
\,
{\rm Op}_\eps^w(
\bd(x,\xi))
\,
U_\eps(t)
\,
\chi_\de\(|x|\)
\Big\|_{{\cal L}(L^2)}
&
\\
&
=
\Big\|
\chi_\de\(|x|\)
{\rm Op}_\eps^w\(
\sum_{k=0}^N \eps^k \bbd^{(k)}(t,x,\xi)
\)
\chi_\de\(|x|\) 
\Big\|_{{\cal L}(L^2)}
&
\\
&
\qquad\qquad\qquad
\qquad\qquad\qquad
\qquad\qquad\qquad
\qquad
+ O_\de\(
\eps^{N+1} \;
\<t\>^{c_\de N^2} \) 
&
\\
&
=
\Big\|
{\rm Op}_\eps^w\(
\chi_\de\(|x|\)\,\sharp\,\(\sum_{k=0}^N \eps^k \bbd^{(k)}(t,x,\xi)
\)
\,\sharp\, \chi_\de\(|x|\)
\)
\Big\|_{{\cal L}(L^2)}
\\
&
\qquad\qquad\qquad
\qquad\qquad\qquad
\qquad\qquad\qquad
\qquad
+
O_\de\(
\eps^{N+1} \;
\<t\>^{c_\de N^2} \) 
.
&
\end{align*}
Now, part (iii) of Lemma \ref{bourob}
and standard symbolic calculus indicate that the above symbol
has support\footnote{
we make here the same abuse of notation than in (\ref{supb}).
}
in
\begin{align*}
&
\mathop{\cup}_{k=0}^N
\({\rm supp}\; \chi_\de\(|x|\)  \cap {\rm supp} \;  \bbd^{(k)}(t,x,\xi)\)
\\
&
\qquad\qquad\qquad
\subset
\{(x,\xi) \text{ s.t. }
|x|\leq \de, \text{ and } 
\Phi(t,x,\xi) \in E(\de) \} .
\end{align*}
The non-trapping condition (and more precisely
estimate (\ref{Xinfty}) below) allows in turn to deduce that this set is
void for $t$ large enough. Hence, up to taking a large value
of $\tau$,
$\tau\geq\tau(\de)$ for some $\tau(\de)$,
we eventually obtain in (\ref{32}),
\begin{align*}
&
\frac{1}{\eps}
\int_{\tau}^{\eps^{-\k}}
\hspace{-0.2cm}
\Big\|
\chi_\de\(|x|\)
\chi_\de\(H_\eps\)
U_\eps(t)
\chi_\de\(|x|\)
\Big\|_{{\cal L}(L^2)}
dt
&
\\
&
\leq
\frac{1}{\eps}
\int_\tau^{\eps^{-\k}}
\hspace{-0.4cm}
O_\de\(
\eps^{(N+1)/2}
\<t\>^{c_\de N^2/2} \) dt
\leq
O_\de\(
\eps^{(N-1)/2- c_\de \k N^2/2 } \)
\leq
O_{\k,\de}(\eps),
&
\end{align*}
for $\k$ small enough (and $N=4$ will do).
This ends the proof of proposition \ref{brob}.


\subsection{Proof of Lemma \ref{bourob}:
an Egorov theorem for polynomially large times}
\label{6.2}


\noi
In view of the above proof,
we are left with the task of proving the large time Egorov
theorem of Lemma \ref{bourob}. To do so, we follow here closely
ideas developped in \cite{BR} in a slightly different context.
Part (iii) of the Lemma is proved in \cite{BR}, so we will skip this aspect.
The implication (i) $\Rightarrow$ (ii) in Lemma \ref{bourob},
which we prove below for completeness, is also essentially proved in \cite{BR}.
Our main task in the sequel turns out to be the proof part (i) of the Lemma.

\medskip

\noi


\noi
The proof is given in several steps.

\medskip

\noi
{\it {\bf First step:} estimates on the flow $\Phi(t,x,\xi)$}

\noi
In this step, we prove that for small enough a $\de$,
there is a time
$T(\de)$, depending on $\de$, such that
for any initial datum
$(x,\xi)$ of phase-space in the set $E(\de)=
\{
|x|\leq \de \}
\cap
\{|\xi^2/2-n^2(x)| \leq \de \}$
(see \ref{supb})),
one has
\bea
\label{Xinfty}
&&
\forall t \geq T(\de) , \quad
|X(t,x,\xi)|
\geq C_{\de} \; t  ,
\eea
for some constant $C_{\de}>0$ that depends on $\de$,
that is however independent
of both time $t$ and the initial point $(x,\xi)$ under consideration.
The proof is standard and uses the information $n_\infty^2>0$.

First, the non-trapping condition implies that for
any large number $R'>0$, and for
any initial point $(x,\xi) \in E(\de)$, there exists a time 
$T(R',x,\xi)$ such that
$$
\forall t \geq
T(R',x,\xi) , \quad
|X(t,x,\xi)|\geq R' .
$$
By continuous dependence of the flow $X(t,x,\xi)$
with respect to the initial data
$(x,\xi)$, and compactness of the set $E(\de)$, there is
a time $T(R',\de)$, that now depends upon $R'$ and $\de$
only, such that
for any initial point $(x,\xi) \in E(\de)$, there holds
$$
\forall t \geq
T(R') , \quad
|X(t,x,\xi)|\geq R' .
$$
In other words, the trajectory $X(t,x,\xi)$ goes to infinity as time
goes to infinity,
uniformly with respect to the initial datum $(x,\xi) \in E(\de)$.

Second, we get estimates for the standard ``escape function'' of
quantum and classical scattering, namely the function $X(t) \cdot \Xi(t)$.
We compute
\begin{align*}
&
\frac{\d}{\d t}
\(X(t,x,\xi) \cdot \Xi(t,x,\xi) \)
=
2\(\frac{\Xi^2(t,x,\xi)}{2}-n^2\(X(t,x,\xi)\)\)
\\
&
\qquad\qquad\qquad\qquad
+2 n^2\(X(t,x,\xi)\)+X(t,x,\xi) \cdot \nabla n^2\(X(t,x,\xi)\)
\\
&
\qquad
=
2 \; \(\frac{\xi^2}{2}-n^2(x)\)
+2 n^2\(X(t,x,\xi)\)+X(t,x,\xi)\cdot \nabla n^2\(X(t,x,\xi)\)
\\
&
\qquad\qquad
\text{ (thanks to the conservation of energy) }
\\
&
\qquad
\mathop{\longrightarrow}\limits_{t \rgt \infty}
2\(\frac{\xi^2}{2}-n^2(x)\)+2 n^2_\infty ,
\end{align*}
uniformly with respect to the initial datum
$(x,\xi) \in E(\de)$. Hence, using the fact that $n^2_\infty>0$,
and taking
a possibly smaller value of the cutoff parameter $\de$,
we obtain the existence
of a constant $C_{\de}>0$, and another time $T(\de)$, such that
\bea
&&
\label{cde}
\forall t \geq T(\de) , \quad
X(t,x,\xi) \cdot \Xi(t,x,\xi)
\geq
C_{\de} \;  t.
\eea
Using the fact that
$
\frac{\d}{\d t}
\(
\frac{1}{2} \; X^2(t,x,\xi)
\)
=
X(t,x,\xi) \cdot \Xi(t,x,\xi) ,
$
we deduce the desired lower bound
$$
\forall t \geq T(\de) , \quad
\frac{1}{2} \( X^2(t,x,\xi) - X^2(T(\de),x,\xi) \)
\geq
C_{\de} \; \frac{t^2}{2} .
$$

\bigskip

\noi
{\it {\bf Second step:} estimates on the linearized flow
$F(t,x,\xi)$.}

\noi
One first proves the estimate (\ref{linflow}) in the case $\a=\b=0$.
By its very definition (\ref{linf}), the linearized flow
$$
F(t,x,\xi)
=
\left(
\begin{array}{cc}
A(t,x,\xi)&B(t,x,\xi)
\\
C(t,x,\xi)&D(t,x,\xi)
\end{array}
\right)
.
$$
satisfies (see (\ref{linf1}), (\ref{linf2})) the differential system
\begin{align}
\non
\label{lin1}
&
\frac{\d}{\d t} A(t,x,\xi)=C(t,x,\xi) ,
&
A(0,x,\xi)={\rm Id} ,
\\
&
\frac{\d}{\d t} C(t,x,\xi)=D^2 n^2\(X(t,x,\xi)\) \; A(t,x,\xi) ,
&
C(0,x,\xi)=0 ,
\end{align}
together with
\begin{align}
\non
\label{lin2}
&
\frac{\d}{\d t} B(t,x,\xi)=D(t,x,\xi) ,
&
B(0,x,\xi)=0 ,
\\
&
\frac{\d}{\d t} D(t,x,\xi)=D^2 n^2\(X(t,x,\xi)\) \; B(t,x,\xi) ,
&
D(0,x,\xi)={\rm Id} .
\end{align}
Here, the notation $D^2 n^2(x)$ refers to the Hessian of the function
$n^2(x)$ in the variable $x$.
Due to the assumption (\ref{decn}) on the behaviour
of $n^2(x)$ at infinity, we readily have
$$
|D^2 n^2(x)| \leq C \; \<x\>^{-\rho-2} ,
$$
for some constant $C>0$, independent of $x$. This, together with the previous
bound (\ref{Xinfty}) on the behaviour of the flow $X(t,x,\xi)$ at infinity
in time, gives the estimate
\bea
\label{decddn}
&&
\big|D^2 n^2\(X(t,x,\xi)\)\big|
\leq
C_0 \;  \<t\>^{-\rho-2} ,
\eea
for some constant $C_0>0$ which is independent of time $t\geq 0$,
and of the point $(x,\xi)$ in phase-space.
We are thus in position to estimate
$A(t)$ and $C(t)$ using (\ref{lin1}).
Integrating (\ref{lin1}) in time, and setting 
\bea
\label{et}
&&
\eps(t):=|D^2 n^2\(X(t,x,\xi)\)|
\eea
for convenience, we obtain (dropping the dependence on $(x,\xi)$
of the various functions),
\begin{align}
\label{esta}
&
|A(t)-{\rm Id}|
\leq
\int_0^t (t-s) \; \eps(s) \; |A(s)-{\rm Id}| \, ds
+
\int_0^t (t-s) \; \eps(s) \, ds
,
\\
\label{estc}
&
|C(t)|
\leq
\int_0^t \eps(s) \; |A(s)| \; ds .
\end{align}
Choose now a constant $C_*$, and define the time $t_*$ as
$$
t_*:=\sup\{t \geq 0 \; \text{ s.t. } \;
|A(t)-{\rm Id}| \leq C_* \<t\>^{1+\eta} \} .
$$
We prove that $t_*=+\infty$, provided $C_*$ is large enough.
Indeed,
for any time $t\leq t_*$, using
(\ref{esta}) together with the decay (\ref{decddn}), we have
\beas
&&
|A(t)-{\rm Id}|
\leq
C_0 C_* \int_0^t (t-s)  \<s\>^{-\rho-1+\eta} \, ds
\leq
C_0 C_* \; t \; \int_0^t \<s\>^{-\rho-1+\eta}\, ds
\\
&&
\qquad
\leq
C_0 C_* C_\eta \; t
\\
&&
\qquad
\;
\text{ (for some constant $C_\eta >0$,
provided $\eta>0$ satisfies $\eta<\rho/2$ ) }
\\
&&
\qquad
<
C_* \; \<t\>^{1+\eta}\\
&&
\qquad \; \text{ (provided
$t$ is large enough, $t\geq T(C_0,C_\eta)$, for some $T(C_0,C_\eta)$}
\\
&&
\qquad\;
\text{ that
only depends on $C_0$ and $C_\eta$).}
\eeas
On the other hand, we certainly have
$|A(t)-{\rm Id}| \leq C_* \<t\>^{1+\eta}$ for bounded values of time
$t\leq T(C_0,C\eta)$, provided $C_*$ is large enough.
Hence $t_*=+\infty$. Inserting this upper-bound for $A$ in (\ref{estc}) gives
$$
|C(t)| \leq C_\eta ,
$$
for some $C_\eta>0$, provided $\eta>0$ is small enough.
We may estimate $B(t)$ and $D(t)$ in the similar way.
The analysis is the same, and starts with the formulae
\beas
&&
|B(t)|
\leq
t  + \int_0^t (t-s) \; \eps(s) \; |B(s)|
\; ds ,
\\
&&
|D(t)|
\leq
1+ \int_0^t \eps(s) \; |B(s)| \; ds .
\eeas
We skip the details.
At this level, we have obtained the bound
$$
|F(t,x,\xi)|
\leq
C_\eta \; \<t\>^{1+\eta} ,
$$
for any (small enough) $\eta>0$, and a constant $C_\eta$ independent of
$(t,x,\xi)$.

\bigskip

\noi
{\it {\bf Third step:}
estimates on the derivatives of the linearized  flow}

\noi
Let now $\a$ be any multi-index.
We prove (\ref{linflow}) by induction on $|\a|$.
Define, for any  $p \geq 1$
$$
M_p(t):=\sup_{|\b|= n} \sup_{(x,\xi)\in \R^{2 d}}
\Bigg|\frac{\d^\b \Phi(t,x,\xi)}{\d (x,\xi)^\b}\Bigg| ,
$$
We have proved in the second step above that
$$
M_1(t)
\leq
C_\eta  \; \<t\>^{1+\eta} .
$$
Assume that for some integer $p_0$, the estimate
$$
M_p(t) \leq C_{p,\eta} \; \<t\>^{p (1+\eta)+2 (p-1)} ,
$$
has been proved for any $p \leq p_0$. We wish to prove the analogous
estimate for $M_{p_0+1}$.
Take any multi-index
$\a$ of length $|\a|=p_0$.
From now on, we systematically omit the dependence of
the various functions and derivatives
with respect to $(x,\xi)$, and write
$\d^\a F(t)$, $\d^\a H$ instead of
$\d^\a F(t,x,\xi)/\d (x,\xi)^\a$,
$\d^\a H(x,\xi)/\d (x,\xi)^\a$
and so on.
Upon differentiating $\a$ times the linearized equation
(\ref{hamlinfs}) on $F$,
we obtain,
\bea
\label{daf}
&&
\d_t \(\d^\a F(t)\)=
J
\sum_{\b\leq\a}
\left(
\begin{array}{c}
\a\\
\b
\end{array}
\right)
\d^\b \(D^2  H\(\Phi(t)\) \) \; \(\d^{\a-\b} F(t)\)
.
\eea
In order to make estimates in (\ref{daf}), we first need
to write the Fa\`a de Bruno formula as
$$
\d^\b \( D^2 H\circ \Phi(t)\)
=
\b !  \; \sum_{\g, m}
\(\d^\g D^2 H \) \circ \Phi(t) \; \times
\prod_\z
\frac{1}{m(\z) !}
\(\frac{\d^\z \Phi(t) }{\z !}\)^{m(\z)} .
$$
Here $\b \in \N^{2d}$,
$\g \in \N^{2d}$, and
$\z \in \N^{2d}$
are multiindices, and $m$ associates to each multi-index
$\z \in \N^{2d}$, another
multi-index $m(\z) \in \N^{2d}$. Also,
the above sum carries over all values of $\g$, $m$, and $\z$ such
that
\bea
\label{cons}
&&
\sum_\z m(\z)= \g , \; \sum_\z \z \; |m(\z)| = \b .
\eea
Finally, when $|\b|\geq 1$, the above sums carries over
$\g$'s and $\z$'s such that $|\g|\geq 1$ and $|\z|\geq 1$.
All this gives in (\ref{daf}),
\beas
&&
\d_t \(\d^\a F(t)\)=
J
\sum_{\b\leq\a}
\b! \;
\left(
\begin{array}{c}
\a\\
\b
\end{array}
\right)
\sum_{\g, m}
\(\d^\g D^2 H\) \circ \Phi(t)
\\
&&
\qquad
\times
\prod_\z
\frac{1}{m(\z) !}
\(\frac{\d^\z \Phi(t)}{\z !}\)^{m(\z)} 
\times \; \d^{\a-\b} F(t)
.
\eeas
Hence, putting apart the contribution stemming from
$\b=0$, we recover
\bea
\label{eqdaf}
&&
\d_t \(\d^\a F(t)\)=
J
\; 
D^2 H\(\Phi(t)\) \; \(\d^\a F(t)\)
+ R_\a(t) ,
\eea
where the remainder term $R_\a(t)$ is estimated by
\begin{align*}
&
|R_\a(t)|
\\
&
\quad
\leq
C_{|\a|} 
\,
\sum_{0\neq\b\leq\a}
\sum_{\g,m}
|\(\d^\g D^2 H\) \circ \Phi(t)|
\prod_\z
\(|\d^\z \Phi(t)|\)^{|m(\z)|}
\, |\d^{\a-\b} F(t)|
&
\\
&
\quad
\leq
C_{|\a|} 
\;
\sum_{0\neq\b\leq\a}
\sum_{\g,m}
\; |\d^{\a-\b} F(t)|
\prod_\z
\(|\d^\z \Phi(t)|\)^{|m(\z)|}
.
&
\end{align*}
for some constant $C_{|\a|}>0$
that depends on $|\a|$. The last line uses the fact that
$\sup_{x,\xi}|\d^\g D^2 H(x,\xi)|\leq C_\g$ for some constant $C_\g$.
Using the inductive assumption,
we recover
\beas
&&
|R_\a(t)|
\leq
C_{|\a|,\eta} \; 
\sum_{0\neq\b\leq\a}
\sum_{\g, m}
\;
\<t\>^{(|\a-\b|+1)(1+\eta)+2|\a-\b|}
\\
&&
\qquad\qquad\qquad\qquad\qquad\qquad
\times
\prod_\z
\<t\>^{\(|\z|(1+\eta)+2(|\z|-1)\) \; |m(\z)| }
\\
&&
\;
\leq
C_{|\a|,\eta} \; 
\sum_{0\neq\b\leq\a}
\<t\>^{ (1+\eta) 
\( 1+|\a-\b|+\sum_\z |\z| |m(\z)| \)
+ 2 \( |\a-\b| + \sum_\z (|\z|-1) \; |m(\z)| \)
}
\\
&&
\;
=
C_{|\a|,\eta} \; 
\sum_{0\neq\b\leq\a}
\<t\>^{(1+\eta) \( 1+|\a-\b|+|\b| \) 
+
2 \( |\a-\b| + |\b|-|\g| \)
 } \;
\\
&&
\leq
C_{|\a|,\eta} \; 
\<t\>^{(1+\eta) \( 1+|\a|\) 
+
2 \( |\a|-1 \)
 } \;
.
\eeas
Here we used the constraints (\ref{cons}) together with the information
$|\g|\geq 1$.
Using Lemma \ref{groj} below in equation (\ref{eqdaf}) satisfied
by $\d^\a F$, we obtain,
$$
|\d^\a F(t) |
\leq
C_{|\a|,\eta} \; 
\<t\>^{(1+\eta) \( |\a|+1 \) + 2 |\a| } .
$$
Hence
$$
M_{p_0+1}(t)
\leq
C_{p_0,\eta} \; 
\<t\>^{(1+\eta) \(p_0 +1 \)+ 2 p_0} .
$$
This ends the recursion.

\bigskip

\noi
{\it{\bf
Fourth step:} A Gronwall
Lemma for solutions to the linearized Hamilton equation}

\noi
The preceding step uses the following
\begin{lemma}
\label{groj}
Assume the function $G(t,x,\xi)$ satisfies the differential equation
\bea
\label{eqg}
&&
\frac{\d G(t,x,\xi)}{\d t} 
= J \cdot D^2 H \(\Phi(t,x,\xi)\) \cdot G(t,x,\xi)
+O\(\<t\>^\l\)
,
\\
&&
\non
G(0,x,\xi)=0 ,
\eea
where the $O\(\<t\>^\l\)$ is uniform in $(x,\xi)$. Then,
$G$ satisfies the uniform estimate
\beas
G(t,x,\xi)=O\(\<t\>^{\l + 2 }\)
.
\eeas
\end{lemma}

\noi
{\bf Proof of Lemma \ref{groj}}

\noi
Decompose $G(t)\equiv G(t,x,\xi)$ as
\beas
G(t)
=
\left(
\begin{array}{cc}
A_G(t) & B_G(t)
\\
C_G(t) & D_G(t)
\end{array}
\right) .
\eeas
Then, equation (\ref{eqg}) for $G$ writes
\begin{align}
\label{lin1g}
&
\frac{\d}{\d t} A_G(t)=C_G(t) + O\(\<t\>^\l\)
,
&
A_G(0)=0 ,
\qquad
\\
&
\non
\frac{\d}{\d t} C_G(t)=D^2 n^2\(X(t)\) \; A_G(t)
+ O\(\<t\>^\l\) ,
&
C_G(0)=0 ,
\qquad
\end{align}
together with
\begin{align}
\label{lin2g}
&
\frac{\d}{\d t} B_G(t)=D_G(t) + O\(\<t\>^\l\)
,
& \quad B_G(0)=0 ,
\qquad
\\
&
\non
\frac{\d}{\d t} D_G(t)=D^2 n^2\(X(t)\) \; B_G(t)
+ O\(\<t\>^\l\)
,
&
D_G(0)=0 .
\qquad
\end{align}
Equations (\ref{lin1g}) give rise to the estimates
\bea
\label{estag}
&&
|A_G(t)|
\leq
C \; \int_0^t 
(t-s) \; \( \; \eps(s) \; |A_G(s)| + \<s\>^\l \) \; ds ,
\\
\label{estcg}
&&
|C_G(t)| \leq C \; \int_0^t \eps(s) \; |A_G(s)| \; ds ,
\eea
where the function $\eps(s)$ is defined in (\ref{et}) above.
Using $\eps(s)\leq C_0 \; \<s\>^{-\rho-2}
\leq C_\eta \; \<s\>^{-\eta-2}$ for any small $\eta>0$
(see (\ref{decddn})),
gives in equation (\ref{estag}),
\bea
\label{estagb}
|A_G(t)|
\leq
C_\eta \; t \; \int_0^t 
\<s\>^{-\eta-2} \; |A_G(s)| \; ds
+
C \; \<t\>^{\l+2} \; 
.
\eea
From this it can be deduced that
$$
|A_G(t)|\leq C \; \<t\>^{\l+2} .
$$
(for a given constant $C_*$, define indeed
$t_*=\sup\{t \geq 0 \; \text{ s.t. } \;
|A_G(t)|\leq C_* \; \<t\>^{\l+2}$ - one deduces from (\ref{estagb})
that $t_*=+\infty$ provided $C_*$ is large enough - see (\ref{esta})
and sequel for details).
Equation (\ref{estcg}) then gives
$$
|C_G(t)|
\leq C_\eta \; \int_0^t \<s\>^{-\eta-2} \; |A_G(s)| \; ds 
\leq C_\eta \; \<t\>^{\l+1-\eta} .
$$
The estimates for $B_G$ and $D_G$ are the same. This ends the proof of
the Lemma.

\bigskip

\noi
{\it {\bf Fifth step}: adapting the estimates of \cite{BR}}

\noi
We now put together the estimates on the linearized flow obtained before,
to complete the proof of parts (ii) and (iii) of Lemma \ref{bourob}.

The construction of the symbols
$\bbd^{(k)}(t,x,\xi)$ in Lemma \ref{bourob}
is made in an explicit way in \cite{BR}.
Part (iii) of Lemma \ref{bourob}
follows. 
Also,
the remainder estimate (\ref{remestbr}) is a consequence
of the above estimates on the linearized
flow $F(t,x,\xi)$ and its derivatives,
upon adapting the analysis of \cite{BR}.
Let us indeed write the rough (but simpler) estimate
$$
|\d^\a F(t,x,\xi)| \leq C_\a \<t\>^{4 |\a| +2} ,
$$
corresponding to the special
choice $\eta=1$ in (\ref{linflow}).
Then, Theorem 1.2 - formula (12) of \cite{BR},
$$
\bbd^{(0)}(t,x,\xi)=
\bd\(\Phi(t,x,\xi)\)
,
$$
together
with the Fa\'a de Bruno formula, give for any multi-index
$\a$ the estimate
$$
|
\d^\a\bbd^{(0)}(t,x,\xi)
|
\leq
C_{|\a|}
\; \<t\>^{4 |\a| } .
$$
From Theorem 1.2 - formula (14) of \cite{BR},
we have for any $k \geq 1$ the explicit value
\begin{align*}
&
\bbd^{(k)}(t,x,\xi)
=
\mathop{\sum}_{
\scriptstyle
\begin{array}{c}
\scriptstyle
|\a|+\ell=k+1\\
\scriptstyle
0\leq \ell \leq k-1
\end{array}
}
\G(\a)
\int_0^t
\[\d^\a H \times \d^\a \bbd^{(\ell)}\]
\circ \Phi(t-s,x,\xi) \, ds
,
&
\end{align*}
where $\G(\a)$ is a harmless coefficient whose explicit value is given in
\cite{BR}. This, together with the Fa\'a de Bruno formula, implies
for any $k \geq 1$, the upper-bound
$$
|
\d^\a\bbd^{(k)}(t,x,\xi)
|
\leq
C_{|\a|,k}
\; \<t\>^{c_0 (k |\a| + k^2+1)},
$$
for some fixed number $c_0$, independent of $\a$ and $k$.
Then, using formulae (51), together with (52), (54), (97) and (99) of
\cite{BR}
gives the estimate (\ref{remestbr}). This ends the proof of
Lemma \ref{bourob}.
\qed


\section{Contribution of moderate times, close to the zero energy}
\label{mod}

\setcounter{equation}{0}


After the work performed in sections \ref{st} through \ref{brbr},
there only remains to estimate the most difficult term
\beas
&&
\frac{1}{\eps}
\int_{\T \eps}^{\tau}
\(1-\chi\)\(\frac{t}{\T \eps}\)
\;
\e^{-\a_\eps t} \;
\Big<
U_\eps(t)
\chi_\de\(H_\eps\) S_\eps
,
\phi_\eps
\Big>
\; dt
.
\eeas
This is the key point of the present paper.

The main result of the present section is the following
\begin{proposition}
\label{dur}
We use the notations of section \ref{redu}.
The index $n^2$ is assumed to have the symbolic behaviour
(\ref{decn}). The zero energy is assumed non-trapping
for the Hamiltonian $\xi^2/2-n^2(x)$.
Finally, we need the tranversality condition {\bf (H)} page \pageref{HH} on the
trajectories $\Phi(t,x,\xi)$ with initial data satisfying
$x=0$, $\xi^2/2=n^2(0)$.
Then, the following two estimates hold true,

(i)
for any {\em fixed} value of the truncation parameters
$\th$, $\tau$ and $\de$,
we have
\beas
&&
\frac{1}{\eps}
\int_{\th}^{\tau}
\(1-\chi\)\(\frac{t}{\th}\)
\;
\e^{-\a_\eps t} \;
\Big<
U_\eps(t)
\chi_\de\(H_\eps\) S_\eps
,
\phi_\eps
\Big>
\; dt
\; 
\mathop{\longrightarrow}\limits_{\eps\rgt 0}
0 .
\eeas

(ii)
for $\th>0$ small enough, there exists a constant $C_\th>0$
such that for any $\eps \leq 1$,
we have
\begin{align*}
\non
&
\frac{1}{\eps}
\int_{\T \eps}^{2\th}
\(1-\chi\)\(\frac{t}{\T \eps}\)
\chi\(\frac{t}{\th}\)
\;
\e^{-\a_\eps t} \;
\Big<
U_\eps(t)
\chi_\de\(H_\eps\) S_\eps
,
\phi_\eps
\Big>
\; dt
\; 
&
\\
&
\qquad\qquad\qquad\qquad
\leq
C_\th \; \T^{-d/2+1} 
\mathop{\longrightarrow}\limits_{\T\rgt +\infty}
0
.
&
\end{align*}
\end{proposition}

\medskip

\noi
The remainder part of this paragraph is devoted to
the proof of proposition \ref{dur}.
In order to shorten the notations, we define
\bea
\wce:=
\(1-\chi\)\(\frac{t}{\T \eps}\)
\;
\e^{-\a_\eps t} ,
\eea
so that the proof of proposition \ref{dur} boils down to estimating
\bea
\label{tob}
&&
\frac{1}{\eps}
\int_{\T \eps}^{\tau}
\wce \,
\Big<
\chi_\de\(H_\eps\) S_\eps
,
U_\eps(-t)
\phi_\eps
\Big>
\; dt
.
\eea
The precise value of the cut-off function $\wce$ in the analysis
of (\ref{tob}) will
be essentially irrelevant in the sequel.

\medskip

\noi
{\bf Proof of proposition \ref{dur}}

The proof is given in several steps.
As in section \ref{brbr}, we begin with some preliminary reductions,
exploiting the informations
on the microlocal support of the various functions.
Then, we use the elegant wave-packet approach of Combescure and Robert
\cite{CRo} to compute the semi-classical propagator $U_\eps(t)$
in (\ref{tob})
in a very explicit way - see Theorem \ref{Crth} below:
this gives a representation in terms of a Fourier
integral operator {\em with complex phase}, that is very well suited for
our asymptotic analysis
(see also \cite{CRR}, or the work by Hagedorn and Joye \cite{H1},
\cite{H2}, \cite{HJ},
or by Robinson \cite{Rb}, or even the seminal work by Hepp \cite{He}
for similar
representations - see also Butler \cite{Bt}).
This eventually reduces the analysis to stationary phase arguments
that are very much in the spirit of \cite{CRR}, and where the tranversality
assumption {\bf (H)} page \pageref{HH} turns out to play a crucial role.

\medskip

\noi
{\it
{\bf First Step:}
Preliminary reduction, projection over the gaussian wave packets}

\noi
As in section \ref{brbr} (see (\ref{32}), (\ref{hchi}), (\ref{supb})),
we may first build up a symbol $a_0(x,\xi) \in C^\infty_c(\R^{2d})$
such that
\bea
\label{a0}
&&
{\rm supp} \, a_0 \subset
\{|x|\leq \de\}\cap \{|\xi^2/2-n^2(x)|\leq  \de\},
\eea
and
\beas
\Big<
\chi_\de\(H_\eps\) S_\eps
,
U_\eps(-t)
\phi_\eps
\Big>
=
\Big<
{\rm Op}^w_\eps\(a_0(x,\xi)\) S_\eps
,
U_\eps(-t)
\phi_\eps
\Big>
+
O_\de(\eps^\infty).
\eeas
With the notation (\ref{hchi}), we actually have the value
$a_0(x,\xi)=\Chid(x,\xi) \sharp \chi_\de\(|x|\)$.
Therefore, the asymptotic analysis of (\ref{tob}) reduces to
that of the expression
\begin{align}
\label{aveca}
&&
\frac{1}{\eps}
\int_{\T \eps}^{\tau}
\wce \,
\Big<
{\rm Op}^w_\eps(a_0)
S_\eps
,
U_\eps(-t)
\phi_\eps
\Big>
\, dt.
\end{align}
Now, to be able to use the wave-packet approach of \cite{CRo}, we need
to decompose the above scalar product on the basis of
the Gaussian wave packets
\beas
\varphi_{q,p}^\eps(x,\xi)
:=
(\pi \eps)^{-d/4}
\; \exp\(\frac{i}{\eps} p \cdot \(x-\frac{q}{2}\)\)
\; \exp\(-\frac{(x-q)^2}{2\eps}\).
\eeas
Each function $\varphi_{q,p}^\eps$ is microlocally
supported near the point $(q,p)$
in phase-space. Using the well-known
orthogonality properties of these states, i.e.
$$
\<u,v\>=(2\pi\eps)^{-d} \int_{\R^{2d}}
dq dp \; \<u,\varphi^\eps_{q,p}\> \; \<\varphi^\eps_{q,p},v\>,
$$
for any $u(x)$ and $v(x)$ in the space $L^2(\R^d)$,
and forgetting the normalizing factors like $\pi$, etc.,
we obtain in (\ref{aveca})
\begin{align}
\non
&
\frac{1}{\eps^{d+1}}
\int_{\T \eps}^{\tau}
\int_{\R^{2d}} dt dq dp \,
\wce
\Big<
{\rm Op}^w_\eps(a_0)
S_\eps
,\varphi^\eps_{q,p}\Big> \,
\Big<\varphi^\eps_{q,p},
U_\eps(-t)
\phi_\eps
\Big>
\\
&
\label{avecac}
=
\frac{1}{\eps^{d+1}}
\int_{\T \eps}^{\tau}
\int_{\R^{2d}} dt dq dp \,
\wce
\Big<
S_\eps
,
{\rm Op}^w_\eps(a_0)
\varphi^\eps_{q,p}\Big> \,
\Big<
U_\eps(t)
\varphi^\eps_{q,p},
\phi_\eps
\Big>
.
\end{align}
Before going further, and in order to prepare for the use of
the stationary phase theorem below,
we make the simple observation that the integral
$dq dp$ over $\R^{2d}$ in (\ref{avecac}) may be carried over the
compact set
$
\{|x|\leq 2\de\}\cap \{|\xi^2/2-n^2(x)|\leq 2 \de\}
$,
up to a negligible error $O_\de(\eps^\infty)$.
For that purpose, take a function
$\ck \in C^\infty_c(\R^{2d})$ such that
\begin{align}
\label{sck}
\non
&
{\rm supp} \, \ck \subset
\{|x|\leq 2\de\}\cap \{|\xi^2/2-n^2(x)|\leq 2\de\}
&
\\
&
\ck
\equiv 1 \; \text{ on } \; 
\left\{
|x|\leq 3\de/2\}\cap \{|\xi^2/2-n^2(x)|\leq 3\de/2
\right\}
.
&
\end{align}
We claim the following estimate holds true:
\bea
\label{gwp}
\int_{\R^{2d}}
dq  dp \,
\(1-\ck\) \,
\Big\|
{\rm Op}^w_\eps(a_0)
\varphi^\eps_{q,p}
\Big\|_{L^2(\R^d)}^2
=O_\de(\eps^\infty) .
\eea
Indeed, we have the following simple computation:
\beas
&&
\big\|
{\rm Op}^w_\eps(a_0)
\varphi^\eps_{q,p}
\big\|_{L^2(\R^d)}^2
=
\big\<
{\rm Op}^w_\eps(a_0\sharp a_0)
\varphi^\eps_{q,p}
,
\varphi^\eps_{q,p}\>
\\
&&
\quad
=
\int_{\R^{2d}}
dx d\xi \,
(a_0 \sharp a_0)(x,\xi) \; W(\varphi^\eps_{q,p})(x,\xi)
\\
&&
\qquad\qquad
\text{ (where $W(\varphi^\eps_{q,p})$ denotes the
Wigner transform of $\varphi^\eps_{q,p}$)}
\\
&&
\quad
=
\eps^{-d}
\int_{\R^{2d}}
dx d\xi \,
(a_0\sharp a_0)(x,\xi) \;  \exp\(-\frac{|q-x|^2+|p-\xi|^2}{\eps}\)
,
\eeas
and the last line uses
the fact that the Wigner transform of $\varphi^\eps_{q,p}$
is a Gaussian. Now, using ${\rm supp}\, (a_0\sharp a_0) \subset
\{|x|\leq \de\} \cap \{|\xi^2/2-n^2(x)|\leq \de\}$, together
with (\ref{sck}), establishes (\ref{gwp}).

Using this estimate (\ref{gwp}),
and replacing back
the factor ${\rm Op}^w_\eps(a_0)$ by the identity in (\ref{avecac}),
we arrive at the conclusion
\begin{align*}
\non
&
\frac{1}{\eps}
\int_{\T \eps}^{\tau}
\wce \,
\Big<
\chi_\de\(H_\eps\) S_\eps
,
U_\eps(-t)
\phi_\eps
\Big>
\; dt
=
O_{\tau,\de}\(\eps^\infty\) +
&
\\
&
\quad
\frac{1}{\eps^{d+1}}
\int_{\T \eps}^{\tau}
\int_{\R^{2d}} dt dq dp \,
\wce \, \ck  \,
\Big<
S_\eps
,
\varphi^\eps_{q,p}\Big> \,
\Big<
U_\eps(t)
\varphi^\eps_{q,p},
\phi_\eps
\Big>
\, .
&
\end{align*}
Our strategy is to now pass to the limit in the term
\bea
\label{avecad}
&&
\frac{1}{\eps^{d+1}}
\int_{\T \eps}^{\tau}
\int_{\R^{2d}} dt dq dp \,
\wce \, \ck  \,
\Big<
S_\eps
,
\varphi^\eps_{q,p}\Big> \,
\Big<
U_\eps(t)
\varphi^\eps_{q,p},
\phi_\eps
\Big>
.
\eea
In order to do so,
we need to compute
the time evolved gaussian wave packet
$U_\eps(t)\varphi^\eps_{q,p}$
in an accurate way.

\bigskip

\noi
{\it {\bf Second Step:}
Computation of $U_\eps(t) \varphi^\eps_{q,p}$ - reducing
the problem
to a stationary phase formula}

\noi
The following theorem is proved in \cite{CRo} (see also \cite{Ro},
\cite{Ro2})
\begin{theorem}{\bf (\cite{CRo}, \cite{Ro})}
\label{Crth}
We use the notations of section \ref{redu}.
Under assumption (\ref{decn}) on the refraction index $n^2(x)$, there
exists a family of functions
$\{p_{k,j}(t,q,p,x)\}_{(k,j)\in \N^2}$, that are
polynomials of degree at most $k$ in the variable $x\in \R^d$,
with coefficients depending on $t$, $q$, $p$,
such that for any
$\eps \leq 1$, the following estimate
holds true: for any given value of $\tau$,
and any given integer $N$, we have, for any time $t\in[0,\tau]$,
\begin{align}
\label{propcr}
\non
&
\Bigg\|
U_\eps(t) \varphi^\eps_{q,p}
-
\exp\(\frac{i}{\eps}\de(t,q,p)\)
{\cal T}_\eps(q_t,p_t) \Lambda_\eps
Q_N(t,q,p,x)
&
\\
&
\quad\quad
{\cal M}(F(t,q,p))
\(\pi^{-d/4} \exp\(-x^2/2\)\)
\Bigg\|_{L^2(\R^d)}
\leq
C_{N,\tau} \; \eps^{N} ,
&
\end{align}
where
\beas
&&
Q_N(t,q,p,x):=
1+\sum_{(k,j)\in I_N} \eps^{\frac{k}{2}-j} p_{k,j}(t,q,p,x)
\, ,
\\
&&
I_N:=\{1\leq j \leq 2 N -1 , \; 1 \leq k-2 j \leq 2 N - 1
, \; k \geq 3 j \}
\, .
\eeas
Here, the following quantities are defined:

\noi
$\bullet$
$\Lambda_\eps$ is the dilation operator
\bea
\label{lam}
\(\Lambda_\eps u\)(x):=\eps^{-d/4} \, u\(\frac{x}{\sqrt{\eps}}\) ,
\eea

\noi
$\bullet$
${\cal T}_\eps(q_t,p_t)$ is the translation (in phase-space) operator
\bea
\label{teps}
\({\cal T}_\eps(q_t,p_t) u\)(x)
:=\exp\(\frac{i}{\eps} p_t \cdot \(x-\frac{q_t}{2}\)\)
\;
u(x-q_t) ,
\eea

\noi
$\bullet$
$(q_t,p_t)$ denotes the trajectory
\bea
\label{qt}
(q_t,p_t):=
\(X(t,q,p),\Xi(t,q,p)\) ,
\eea

\noi
$\bullet$
$\de(t,q,p)$ denotes quantity
\bea
\label{delt}
\de(t,q,p)=
\int_0^t \(\frac{p_s^2}{2}+n^2(q_s) \) \; ds
-\frac{q_t \cdot p_t - q \cdot p}{2} ,
\eea

\noi
$\bullet$
${\cal M}(F(t,q,p))$ is the metaplectic operator associated with the
symplectic matrix
$F(t,q,p)$. It acts on the gaussian as
\begin{align}
\label{metf}
&
\non
{\cal M}(F(t,q,p))
\(\exp\(-\frac{x^2}{2}\)\)
=
&
\\
&
\qquad\qquad
\det(A(t,q,p)+iB(t,q,p))^{-1/2}_{\rm c}
\;
\exp\(i \frac{\G(t,q,p) x \cdot x}{2} \)
.
&
\end{align}
Here, the square root
$\det(A(t,q,p)+iB(t,q,p))^{-1/2}_{\rm c}$ is defined
by continuously (hence the index ${\rm c}$)
following the argument of the complex number
$\det(A(t,q,p)+iB(t,q,p))$ starting from its value $1$ at time $t=0$.
Also, the complex matrix $\G(t,q,p)$ is defined as
\bea
\label{gt}
\G(t,q,p)
=
(C(t,q,p)+iD(t,q,p))
\; (A(t,q,p)+iB(t,q,p))^{-1} .
\eea
\end{theorem}

\noi
{\bf Remark}

\noi
If the refraction index $n^2(x)$ is {\em quadratic} in $x$,
then formula (\ref{propcr}) is {\em exact}, and the
whole family $\{p_{k,j}\}$ vanishes. This is essentially
a consequence of the Mehler formula. We refer to \cite{Fo}
for a very complete discussion about the propagators
of pseudo-differential operators with {\em quadratic} symbols.

In the case when $n^2(x)$ is a general function, the polynomials
$p_{k,j}$ are obtained in \cite{CRo}
using perturbative expansions ``around the quadratic case''.
We refer to \cite{Ro} for a very clear and elegant
derivation of these polynomials. Let us quote that similar
formulae are derived and used in \cite{HJ}. The idea
of considering such perturbations ``around the quadratic case'' traces
back to \cite{He}, see also \cite{H1}, \cite{H2}, \cite{Rb}.

The fact that the matrix $A(t)+iB(t)$ is invertible, and
$\G(t)$ is well defined, is proved in \cite{Fo}, see also \cite{Ro2}.
It is a consequence
of the symplecticity of $F(t)$ (see the relations (\ref{relsym})).
We refer to the sequel for an explicit use
of these important relations.
\qed

\medskip

\noi
In the next lines, we apply the above theorem,
and transform formula (\ref{avecad}) accordingly.

On the one hand,
we use the Parseval formula in (\ref{avecad})
to compute
the two scalar products. Forgetting
the normalizing factors like $\pi$, etc., it gives, e.g. for the first
scalar product,
\beas
&&
\<S_\eps,\varphi^\eps_{q,p}\>
=
\eps^{-d/2}
\int_{\R^d} dx d\xi
\,
\exp(i x\cdot \xi/\eps)
\,
{\wh S}(\xi)
\,
\varphi^\eps_{q,p}(x)
\\
&&
\qquad
=
\eps^{-d/2}
\int_{\R^d} dx d\xi
\,
\exp(i x\cdot \xi/\eps)
\,
\chi_1(x)
\,
{\wh S}(\xi)
\,
\varphi^\eps_{q,p}(x)
+O(\eps^\infty) ,
\eeas
for any truncation function $\chi_1$ being $\equiv 1$ close
to the origin.
On the other hand,
we use formula (\ref{propcr})
to compute $U_\eps(t) \varphi^\eps_{q,p}$ in (\ref{avecad}),
using the short-hand notation
\beas
P_N(t,q,p,x):=
\pi^{-d/4}
\det(A(t,q,p)+iB(t,q,p))^{-1/2}_{\rm c}
\,
Q_N(t,q,p,x)
.
\eeas

These two tasks being done, we eventually obtain
in (\ref{avecad}), upon
computing the relevant phase factors explicitly,
\begin{align}
\label{tobd}
\non
&
\frac{1}{\eps}
\int_{\T \eps}^{\tau}
\wce \,
\Big<
\chi_\de\(H_\eps\) S_\eps
,
U_\eps(-t)
\phi_\eps
\Big>
\; dt
=
O_{\tau,\de}\(\eps^{\infty}\) 
+
&
\\
\non
&
\quad
\frac{1}{\eps^{(5d+2)/2}}
\int_{\T \eps}^{\tau}
\int_{\R^{6d}}
dt dq dp d\xi  d\eta  dx dy \,
\wce
\exp\(\frac{i}{\eps}
\Psi(x,y,\xi,\eta,q,p,t)\)
&
\\
&
\qquad\qquad
\wh S(\xi)
\wh\phi^*(\eta)
\ckk P_N\(t,q,p,\frac{y-q_t}{\sqrt\eps}\)
.
&
\end{align}
where $\chi_1\in C^\infty_c$ is $\equiv 1$ close to $(0,0)$.
Here, the crucial (complex) phase factor has the value
\begin{align}
\label{close}
\non
&
\Psi(x,y,\xi,\eta,q,p,t)
=
\int_0^t \(\frac{p_s^2}{2}+n^2(q_s) \) \, ds
-p\cdot\(x-q\)
+p_t\cdot\(y-q_t\)
&
\\
&
\qquad
+x\cdot\xi
-y\cdot\eta
+i\frac{(x-q)^2}{2}
+\frac{\G(t) (y-q_t) \cdot (y-q_t)}{2}
&
\end{align}
Our goal is now to apply the stationary phase formula
to estimate (\ref{close}).
Obviously, the cutoff in time away from $t=0$ in
(\ref{tobd}) prevents one to use directly the stationary
phase formula close to $t=0$. This is the reason
why times close to $0$ are treated apart in the sequel
(see steps four and five below - see also the outline of proof
given in section \ref{redu}).

\bigskip

\noi
{\it {\bf Third Step:}
computing the first and second order derivatives of the phase
$\Psi$}

\noi
First, it is an easy exercice, using the symplecticity relations
(\ref{relsym}), to prove that the matrix $\G(t)$ is symmetric
and it has positive imaginary part. The relation
$$
{\rm Im}\,
\(\G(t) (y-q_t) \cdot (y-q_t)\)
=
\big|
\(A(t)+iB(t)\)^{-1} \, (y-q_t)
\big|^2 ,
$$
implies indeed
\beas
&&
{\rm Im} \, \Psi
=
|x-q|^2+\big|
\(A(t)+iB(t)\)^{-1} \,(y-q_t)
\big|^2 .
\eeas
Hence we recover the equivalence
\bea
&&
\label{im}
{\rm Im}\,\Psi=0 \; \text{ iff } \;
y=q_t \; \text{ and } \; x=q .
\eea
Second,
using the differential
system (\ref{linf1}), (\ref{linf2}) satisfied by the matrices
$A(t)$, $B(t)$, $C(t)$, and $D(t)$, we prove
\beas
&&
\nabla_{q,p}
\(
\int_0^t \(\frac{p_s^2}{2}+n^2(q_s) \) \, ds
\)
=
\left(
\begin{array}{c}
A(t)^\tr p_t-p
\\
B(t)^\tr p_t
\end{array}
\right)
.
\eeas
This gives the value of the gradient of $\Psi$
\bea
\label{nablapsi}
&&
\non
\nabla_{x,y,\xi,\eta,q,p,t}\Psi(x,y,\xi,\eta,q,p,t)
=
\\
&&
\quad
\left(
\begin{array}{c}
-p+\xi+i(x-q)
\\
p_t-\eta+\G(t)(y-q_t)
\\
x
\\
-y
\\
C(t)^\tr (y-q_t)
+ i (q-x)
+A(t)^\tr \G(t) (q_t-y)
\\
-(x-q)
+D(t)^\tr (y-q_t)
+B(t)^\tr \G(t) (q_t-y)
\\
-\disp\frac{p_t^2}{2}+n^2(q_t)
+ \nabla n^2(q_t) \cdot (y-q_t)
+  p_t \cdot \G(t) (q_t-y)
\end{array}
\right)
.
\eea
This computation is done
up to irrelevant $O\( (y-q_t)^2+(x-q)^2 \)$ terms.

\medskip

\noi
These observations
allow to compute the stationary set, defined as
\begin{align}
\non
&
M:=
\{(x,y,\xi,\eta,q,p,t) \in \R^{6d}\times ]0,+\infty[
&
\\
&
\qquad
\qquad
\; \text{s.t. } \;
{\rm Im}\,\Psi=0  \; \text{ and } \; \nabla_{x,y,\xi,\eta,q,p}\Psi=0
\} .
&
\end{align}
Note (see above) that we exclude the original time $t=0$
in the definition of
$M$. In view of (\ref{im}), (\ref{nablapsi}),
the set $M$ has the value
\bea
\label{crit} 
\non
&&
\non
M=
\left\{
(x,y,\xi,q) \; \text{ s.t. } \;
x=y=q=0
, \;
\xi=p
\right\}
\\
&&
\qquad\qquad
\cap
\;
\left\{
(p,\eta,t) \; \text{ s.t. } \;
\frac{\eta^2}{2}=n^2(0)
, \;
q_t=0 , \; p_t=\eta
\right\}
.
\eea
Note that the second set reads
also, by definition,
$$
\left\{
(p,\eta, t) \; \text{ s.t. } \;
\frac{\eta^2}{2}=n^2(0) , \;
X(t,0,p)=0 , \; \Xi(t,0,p)=\eta
\right\} .
$$

\medskip
\noi
Last, there remains to compute the Hessian of $\Psi$ at the stationary
points. A simple but tedious computation gives, for any
point $(x,y,\xi,\eta,q,p,t)\in M$, the value
\nopagebreak
\newcommand{\scr}{\scriptstyle}
\begin{align*}
\non
&
D^2_{x,y,\xi,\eta,q,p,t}\Psi\,\Big|_{(x,y,\xi,\eta,q,p,t)\in M}
=
\end{align*}
\begin{align*}
\non
&
\left(
\begin{array}{lllllll}
\scr
i {\rm Id}
&\scr
0
&\scr
{\rm Id}
&\scr
0
&\scr
-i {\rm Id}
&\scr
-{\rm Id}
&\scr
0
\\
\scr
0
&\scr
\G_t
&\scr
0
&\scr
-{\rm Id}
&\scr
C_t
-\G_t A_t
&\scr
D_t-\G_t B_t
&\scr
\nabla n^2(0)
\\
\scr
&\scr
&\scr
&\scr
&\scr
&\scr
&\scr
\;\;\;
-\G_t \eta 
\\
\scr
{\rm Id}
&\scr
0
&\scr
0
&\scr
0
&\scr
0
&\scr
0
&\scr
0
\\
\scr
0
&\scr
-{\rm Id}
&\scr
0
&\scr
0
&\scr
0
&\scr
0
&\scr
0
\\
\scr
-i{\rm Id}
&
\scr
C_t^\tr-A_t^\tr \G_t
&\scr
0
&\scr
0
&\scr
-C_t^\tr A_t
+i {\rm Id}
&\scr
-C_t^\tr B_t
&\scr
-C_t^\tr \eta
\\
&\scr
&\scr
&\scr
&\scr
\;\; \;\;
+A_t^\tr \G_t A_t
&\scr
\;\; \;
+A_t^\tr \G_t B_t
&\scr
\;\; \;
+A_t^\tr\G_t \eta
\\
\scr
-{\rm Id}
&\scr
D_t^\tr -B_t^\tr\G_t
&\scr
0
&\scr
0
&\scr
{\rm Id}
-D_t^\tr A_t
&\scr
-D_t^\tr B_t
&\scr
-D_t^\tr \eta
\\
\scr
&\scr
&\scr
&\scr
&\scr
\;\; \;
+B_t^\tr \G_t
A_t
&\scr
\;\;\;
+B_t^\tr \G_t B_t
&\scr
\;\;\;
+B_t^\tr \G_t \eta
\\
\scr
0
&\scr
\nabla n^2(0)^\tr
&\scr
0
&\scr
0
&\scr
-\eta^\tr C_t
&\scr
-\eta^\tr D_t
&\scr
-\eta^\tr \nabla n^2(0)
\\
\scr
&\scr
\; \; \;
- \eta^\tr \G_t 
&\scr
&\scr
&\scr
\;\; \;
+\eta^\tr \G_t A_t
&\scr
\;\; \;
+\eta^\tr \G_t B_t
&\scr
\;\;\;
+\eta^\tr \G_t \eta
\end{array}
\right).
&
\\
&&
\end{align*}
Here we wrote systematically $A_t$, $B_t$, etc. instead of
$A(t)$, $B(t)$, etc. The above matrix is symmetric,
due to the relation (\ref{uv}).
The very last computation we need is that
of ${\rm Ker}\, D^2\Psi$ at stationary points. The value
of $D^2 \Psi|_M$ clearly shows that
\beas
\non
{\rm Ker}\, \(D^2\Psi|_M\)
&=&
\Big\{
(X,Y,\Xi,H,Q,P,T)
\; \text{ s.t. } \; 
X=Y=Q=0 , \; \Xi=P ,
\\
&&
\qquad
-H+(D_t-\G_t B_t) P +T (\nabla n^2(0)-\G_t \eta)=0 ,
\\
&&
\qquad
(-C_t^\tr +A_t^\tr \G_t) B_t P
+T (-C_t^\tr+A_t^\tr \G_t) \eta =0 ,
\\
&&
\qquad
(-D_t^\tr+B_t^\tr \G_t) B_t P
+T (-D_t^\tr+B_t^\tr \G_t)\eta=0 ,
\\
&&
\qquad
\eta^\tr(-D_t+ \G_t B_t)P
+T \eta^\tr
(-\nabla n^2(0)+ \G_t \eta)
=0 \Big\}
.
\eeas
Hence, using
$D_t^\tr-B_t^\tr \G_t=(A_t+i B_t)^{-1}$, together with
$C_t^\tr-A_t^\tr \G_t=-i(A_t+i B_t)^{-1}$,
and $(A_t+iB_t)^{-1,\tr}+\G_t B_t=D_t$
(see (\ref{uv})), we obtain
\bea
&&
\non
{\rm Ker}\, \(D^2\Psi|_M\)
=
\Big\{
(X,Y,\Xi,H,Q,P,T)
\; \text{ s.t. } \; 
X=Y=Q=0 , \; \Xi=P ,
\\
&&
\; \text{ and }
\eta^\tr H = 0 , \;
B_t P + T \eta =0 ,
H=D_t P +T \nabla n^2(0)=0
\Big\}
.
\eea

\bigskip

\noi
{\it{\bf Fourth Step:}
Application of the stationary phase Theorem
- proof of part (i) of proposition \ref{dur}}

\noi
In this step, we formulate the main geometric
assumption on the flow $\Phi(t,x,\xi)$, that allows for the proof
that the contribution in (\ref{tobd}) vanishes asymptotically.

\label{HH}
\medskip
{\centerline
{\bf 
\underline{
(H) Transversality assumption on the flow}
}
}
\medskip
{\it

We suppose that the stationary set
\beas
M=
\left\{
x=y=q=0
, \,
\xi=p
\right\}
\cap
\left\{
\frac{\eta^2}{2}=n^2(0)
, \,
X(t,0,p)=0, \, \Xi(t,0,p)=\eta
\right\}
\eeas

is a {\em smooth submanifold}
of $\R^{6d}\times ]0,+\infty[$, satisfying the additional constraint
\bea
\label{93}
&&
k:={\rm codim} M>5d+2 .
\eea

We also assume that at each point
$m=(x,y,\xi,\eta,q,p,t)\in M$,

the tangent space of $M$ at $m$ is
\bea
\label{94}
\non
&&
{\rm T}_m M
=
\Big\{
(X,Y,\Xi,H,Q,P,T)
\; \text{ s.t. } \; 
X=Y=Q=0 , \; \Xi=P ,
\\
&&
\;
\text{ and }
\eta^\tr H = 0 , \;
B_t P + T \eta =0 ,
-H+D_t P +T \nabla n^2(0)=0
\Big\}
.
\eea

In other words, we assume that ${\rm T}_m M$ is precisely given by

linearizing the equations defining $M$. 
}

\medskip
\noi
{\bf First Remark}

\noi
We show below examples of flows satisfying the above assumption.
It is a natural, and generic, assumption.
Note in particular that the assumption on the codimension
is natural, in that the equations defining $M$ give (roughly)
$4d$ constraints on $(x,y,q,\xi)$, one constraint on $\eta$,
and again $2d$ constraints on the momentum
$p$, the solid angle $\eta/|\eta|$, and time $t$. Hence
one has typically $k=6d+1$.
\qed
\medskip

\noi
{\bf Second remark}

\noi
Equivalently, the above assumption
may be formulated as follows.
The set 
$$
{\cal M}:=
\{(p,\eta,t) \text{ s.t. } \frac{\eta^2}{2}=n^2(0) , \;
X(t,0,p)=0 , \; \Xi(t,0,p)=\eta\}
$$
is assumed to be a smooth
submanifold of $\R^{2d+1}$, satisfying the additional constraint
$
{\rm codim} \, {\cal M}
>
d+2,
$
and whose tangent space is given by
$$
\{(P,H,T) \text{ s.t. }
\eta^\tr H=0 , \;
B_t P + T \eta = 0 , \; 
D_t P + T \nabla n^2(0)-H=0\}.
$$
Note in passing that the conservation of energy allows
to replace the requirement $\eta^2/2=n^2(0)$ by the equivalent
$p^2/2=n^2(0)$ in the definition of ${\cal M}$.
\qed
\medskip

\noi
{\bf Third Remark}
\nopagebreak

\noi
Provided $M$ is a smooth submanifold
with tangent space given upon linearizing the constraints, its codimension
anyhow satisfies
$$
{\rm codim} \, M \geq 5d+2.
$$
Equivalently, provided ${\cal M}$ is a smooth submanifold with
the natural tangent space, its codimension anyhow satisfies
$$
{\rm codim}\,{\cal M} \geq d+2.
$$
As a consequence, the analysis given below (see (\ref{stp}))
establishes that $\<w^\eps,\phi\>$ is uniformly bounded
in $\eps$.
This fact is not known in the literature.
\qed

\medskip
\noi
Under assumption {\bf (H)}, we are ready to use the stationary phase
Theorem in (\ref{tobd}), at least for large enough times  $t$
(recall that the
very point $t=0$ is excluded from the definition of $M$ above).
Indeed, assumption {\bf (H)} precisely asserts the equality
$$
{\rm T}_m M={\rm Ker}\, \(D^2\Psi|_M\) ,
$$
so that the Hessian $D^2\Psi|_M$ is non-degenerate on the normal
space $({\rm T}_m M)^\perp$. This is exactly the non-degeneracy
that we need in order to apply the stationnary phase Theorem.

To perform the claimed stationary phase argument, we first
take a (small) parameter
$$
\th >0 .
$$
We use a cutoff in time $\chi(t/\th)$ with $\chi$
as in (\ref{chi}), and evaluate the contribution
\begin{align*}
&
\frac{1}{\eps}
\int_{\th}^{\tau}
\wce \, \(1-\chi\(\frac{t}{\th}\)\)
\Big<
\chi_\de\(H_\eps\) S_\eps
,
U_\eps(-t)
\phi_\eps
\Big>
\; dt
=
O_{\tau,\de}\(\eps^{\infty}\) +
&
\\
\non
&
\quad
\frac{1}{\eps^{(5d+2)/2}}
\int_{\th}^{\tau}
\int_{\R^{6d}}
\wce \, \(1-\chi\(\frac{t}{\th}\)\)\,
\exp\(\frac{i}{\eps}
\Psi(x,y,\xi,\eta,q,p,t)\)
\\
&
\qquad\quad
\wh S(\xi)
\wh\phi^*(\eta)
\ckk P_N\(t,q,p,\frac{y-q_t}{\sqrt\eps}\)
\, dt dx dy d\xi d\eta dq dp
.
&
\end{align*}
When the point $(x,y,\xi,\eta,q,p,t)$ is far from the stationary
set $M$, the integral is $O(\eps^\infty)$. Close
to the stationary
set $M$, using the fact that the integral carries over a compact support,
we may use
a partition of unity close to $M$, and on each piece
we may use straightened coordinates $(\a,\b)\in\R^{6d+1-k}\times\R^k$
such that
\begin{align*}
&
(x,y,\xi,\eta,q,p,t)=\g(\a,\b) , \;
\text{ where } \; 
\g \; \text{ is a local diffeomorphism, with }
&
\\
&
\qquad
(x,y,\xi,\eta,q,p,t)\in M \Longleftrightarrow \a=0 .
&
\end{align*}
Using such coordinates, we recover a finite sum of terms of the form
\begin{align}
\non
&
\frac{1}{\eps^{(5d+2)/2}}
\int_{\Omega}
dx dy d\xi d\eta dq dp
\,
\exp\(\frac{i}{\eps}
\Psi(x,y,\xi,\eta,q,p,t)\)
\\
\non
&
\quad\quad
\wh S(\xi)
\wh\phi^*(\eta)
P_N\(t,q,p,\frac{y-q_t}{\sqrt\eps}\)
\chi_2(x,y,\xi,\eta,q,p,t)
&
\\
\non
&
=
\frac{1}{\eps^{(5d+2)/2}}
\int_{\Omega'\times\Om''}
d\a\,d\b\,
\exp\(\frac{i}{\eps}
\Psi\circ\g(\a,\b)\)
&
\\
\label{exp1}
&
\qquad\qquad
\(
\wh S(.)
\wh\phi^*(.)
P_N\(.,.,.,\frac{.}{\sqrt\eps}\)
\)\circ\g(\a,\b) \;
\chi_3(\a,\b) ,
&
\end{align}
where $\Om$, $\Om'$, $\Om''$
are bounded, open subsets, and $\chi_2$, $\chi_3$
are cutoff functions. Thanks to the non-degeneracy of the Hessian
$D^2\Psi$ in the normal direction to $M$, for any $\b$,
we have
$$
\(
{\rm det}
\frac{D^2 \Psi\circ \g}{D \a^2}
\)(0,\b) \neq 0 .
$$
Hence,
by the standard stationary phase Theorem,
for any integer $J$, the above integral has the asymptotic expansion
to order $J$
\begin{align}
\non
&
\hspace{-0.5cm}
\eps^{(k-5d-2)/2}
\int_{\Omega''}
d\b\,
\exp\(\frac{i}{\eps}
\Psi\circ\g(0,\b)\)
\\
&
\label{exp2}
\times
\sum_{j=0}^J
\eps^j\;
Q_{2j}(\d_\a,\d_\b)
\(
\(
\wh S(.)
\wh\phi^*(.)
P_N\(.,.,.,\frac{.}{\sqrt\eps}\)
\)\circ\g \;
\chi_3
\)
(0,\b)
\\
&
\non
\hspace{-0.8cm}
+
\eps^{(k-5d-2)/2}
O\Bigg(
\eps^{J+1}
\mathop{\sup}_{k \leq 2 J +d+3}
\Big\|
\d^k_{(\a,\b)}
\Big(
\wh S(.)
\wh\phi^*(.)
P_N\(.,.,.,\frac{.}{\sqrt\eps}\)
\chi_3
\Big)
\Big\|
\Bigg)
,
\end{align}
where the $Q_{2j}$'s are differential operators of order $2j$.
Now, we anyhow have
\beas
\forall
j \in \N \, \qquad
\eps^j \d_y^{2j} P_N\(.,.,.,\frac{y}{\sqrt\eps}\)=O(1) .
\eeas
On the more, $P_N$ is a {\em polynomial} of degree $\leq 4N$
in its last argument. This implies that
the $\eps^{(k-5d-2)/2}O(\ldots)$ in (\ref{exp2}) has at most the size
$$
O\(\eps^{J+1+(k-5d-2)/2-2 N}\) .
$$
Hence, taking $J$ large enough ($J\geq 2 N$ will do),
we eventually obtain in (\ref{exp2}), using
the assumption {\bf (H)} on the codimension $k$ ($k>5d+2$),
\bea
\label{stp}
&&
\non
\frac{1}{\eps^{(5d+2)/2}}
\int_{\th}^{\tau}
\wce \, \(1-\chi\(\frac{t}{\th}\)\)
\Big<
\chi_\de\(H_\eps\) S_\eps
,
U_\eps(-t)
\phi_\eps
\Big>
\; dt
\\
&&
\qquad\qquad
=
O_{\th,\tau,\de}\(\eps^{(k-5d-2)/2}\)
\mathop{\longrightarrow}\limits_{\eps\rgt 0} 0
.
\eea

\bigskip
\noi
{\it{\bf Fifth Step:}
Elimination of times such that  $\T \eps \leq t \leq \th$
- proof of part (ii) of proposition \ref{dur}}

\noi
The previous step leaves us with the task of estimating
\begin{align}
\non
&
\frac{1}{\eps}
\int_{\T \eps}^{2\th}
\wce \, \chi\(\frac{t}{\th}\)
\Big<
\chi_\de\(H_\eps\) S_\eps
,
U_\eps(-t)
\phi_\eps
\Big>
\; dt
.
\end{align}
The idea is to now come back
to the semiclassical scale, and write
\begin{align}
\non
&
\label{cont}
\frac{1}{\eps}
\int_{\T \eps}^{2\th}
\wce \, \chi\(\frac{t}{\th}\)
\Big<
\chi_\de\(H_\eps\) S_\eps
,
U_\eps(-t)
\phi_\eps
\Big>
\; dt
&
\\
&
\;
=
\int_{\T}^{2\th/\eps}
\chi\(\frac{\eps t}{\th}\)
\Big<
\chi_\de\(H_\eps\) S_\eps
,
\exp\(-i t \(\eps^2\D+n^2(x)\)\)
\phi_\eps
\Big>
\, dt
\, .
&
\end{align}
This term is expected to be small, provided $\T$ is large
enough. Indeed, the propagator $\exp\(-i t \(\eps^2\D+n^2(x)\)\)$
acting on $\phi_\eps$
is expected to be close to the free propagator
$\exp\(-i t \(\eps^2\D+n^2(0)\)\)$ on the time-scale we consider.
Hence the propagator should
have size $O(t^{-d/2})$ for large values of time,
and the above time integral
should be $O(\T^{-d/2+1})\rgt 0$ as $\T \rgt \infty$.

We give below a quantitative proof of this rough statement, based
on the exact computation of the propagator $\exp\(-i t \(\eps^2\D+n^2(x)\)\)$
obtained in Theo\-rem \ref{Crth}. The proof given below
could easily be replaced by a slightly simpler one, upon
writing the propagator as a Fourier Integral Operator with {\em real} phase.
We do not detail this aspect, since we anyhow had to use in
the previous steps the more precise
expansion of the propagator given by Theorem \ref{Crth}: this theorem
has indeed the great advantage to give a representation
of the propagator that
is valid {\em for all times}.

From the second step above (see (\ref{tobd})), we know
\begin{align}
\label{zz}
\non
&
\int_{\T}^{2\th/\eps}
\chi\(\frac{\eps t}{\th}\)
\Big<
\chi_\de\(H_\eps\) S_\eps
,
\exp\(-i t \(\eps^2\D+n^2(x)\)\)
\phi_\eps
\Big>
\; dt
&
\\
\non
&
=
O_{\tau,\de}\(\eps^\infty\)
+
\int_{\T}^{2\th/\eps}
\chi\(\frac{t}{\th}\)\,
\times
\eps^{-\frac{5d}{2}}
\int_{\R^{6d}}
\exp\(i\Psi(\eps t)/\eps\)
&
\\
&
\;
\wh S(\xi)
\wh\phi^*(\eta)
\ckk P_N\(t,q,p,\frac{y-q_{\eps t}}{\sqrt\eps}\)
\, dx dy d\xi d\eta dq dp
,
&
\end{align}
where we drop the dependence of the phase
$\Psi$ in $(x,y,\xi,\eta,q,p)$. To estimate this term,
we now concentrate our attention on the space integral
\begin{align}
\label{ft}
&
\non
f_\eps(t):=
\eps^{-\frac{5d}{2}}
\int_{\R^{6d}}
\exp\(i\frac{\Psi(\eps t)}{\eps}\)
\;
\wh S(\xi)
\wh\phi^*(\eta)
&
\\
&
\qquad\qquad
\ckk P_N\(t,q,p,\frac{y-q_{\eps t}}{\sqrt\eps}\)
\, dx dy d\xi d\eta dq dp
.
&
\end{align}
We claim we have the following dispersion estimate,
{\em uniformly in} $\eps$,
\bea
\label{clai}
\big|
f_\eps(t)
\big|
\leq
C_\th \; t^{-d/2} ,
\; \; 
\text{ for some $C_\th>0$, provided }
\;
\T \leq t \leq 2\th/\eps .
\eea
Assuming (\ref{clai}) is proved, equation (\ref{zz}) shows that
\begin{align}
&
\disp
\frac{1}{\eps}
\Bigg|
\int_{\T \eps}^{2\th}
\wce \, \chi\(\frac{t}{\th}\)
\Big<
\chi_\de\(H_\eps\) S_\eps
,
U_\eps(-t)
\phi_\eps
\Big>
\, dt
\Bigg|
\leq
C_\th \, \T^{-\frac{d}{2}+1}
\mathop{\longrightarrow}\limits_{\T\rgt \infty}0 \, ,
&
\end{align}
in any dimension $d\geq 3$, which is enough for our purposes.
It is thus sufficient to prove (\ref{clai}).

\medskip

We have in mind that the integral (\ref{ft})
defining $f_\eps(t)$
should concentrate on the set $x=y=q=0$, $q_t=0$,
$p_t=\eta$, $p=\xi$. Also, the present case should be close
to the ``free'' case where the refraction index $n^2(x)$ has frozen
coefficients at the origin $n^2(x)\approx n^2(0)$.
For that reason, we perform in (\ref{ft})
the changes of variables
\beas
&&
(x-q)/\sqrt\eps \rgt x , \;
(y-q_{\eps t})/\sqrt\eps \rgt y , \;
q \rgt \sqrt\eps q , \;
\\
&&
\xi \rgt p+\sqrt\eps \xi , \;
\eta \rgt \Xi(\eps t,\sqrt\eps q,p)+\sqrt\eps \eta
.
\eeas
We also put apart the important phase factors in the
obtained formula. This gives
\begin{align}
\label{ftt}
f_\eps(t)
=
\int_{\R^{4d}}
dq dp d\eta\,
\exp\Big(
i t \; \wt{\Psi}(p,\eps t,\sqrt{\eps}q,\sqrt{\eps}\eta)
\Big)
\;
G(q,p,\eta,\eps t,\sqrt{\eps}q,\sqrt{\eps}\eta) ,
\end{align}
up to introducing the phase
\begin{align*}
&
\wt{\Psi}(p,\eps t,\sqrt{\eps}q,\sqrt{\eps}\eta)
:=
\frac{1}{\eps t}
\int_0^{\eps t}
\(
\frac{\Xi(s,\sqrt\eps q,p)^2}{2}
+n^2\(X(s,\sqrt\eps q,p)\)
\)
ds
&
\\
\non
&
\qquad
+
\frac{\sqrt\eps p \cdot q
-
\Xi(\eps t,\sqrt\eps q,p)\cdot X(\eps t,\sqrt\eps q,p)
}{\eps t}
&
\\
&
\qquad
+
\sqrt{\eps}\eta
\cdot
\frac{
\sqrt{\eps} q
-
X(\eps t,\sqrt\eps q,p)
}
{\eps t}
\, ,
&
\end{align*}
together with the  amplitude ($C^\infty$, and compactly supported in $p$,
$\sqrt{\eps}q$)
\begin{align}
\label{bigg}
\non
&
G(q,p,\eta,\eps t,\sqrt{\eps}q,\sqrt{\eps}\eta)
:=
\int_{\R^{3d}}
dx dy d\xi\,
\exp\(i\xi\cdot(q+x)-i \eta\cdot (y+q)\)
&
\\
\non
&
\qquad
\exp\(
-\frac{x^2}{2}
+i\frac{\G(\eps t,\sqrt\eps q,p) y \cdot y}{2}
\)
&
\\
\non
&
\qquad
\wh S(p+\sqrt\eps\xi)
\wh\phi^*(\Xi(\eps t,\sqrt\eps q,p)+\sqrt\eps\eta)
\chi_0\(\sqrt{\eps}q,p\)
\\
&
\qquad
\chi_1(\sqrt{\eps}(q+x),X(\eps t,\sqrt\eps q,p)+\sqrt\eps y)
\;
P_N\(t,\sqrt\eps q,p,y\)
.
&
\end{align}
Now, the  idea is to use the stationary phase formula in the $p$
variable in (\ref{ftt}),
where $t$ plays the role of the large parameter.
We wish indeed to recognize in (\ref{ftt}) a formula of the form
$$
\int dp \, \exp\(-it\frac{p^2}{2}\) \times \text{smooth}(p)
,
$$
to recover the claimed decaying factor $t^{-d/2}$ in (\ref{clai}).
In other words, we wish to get the same dispersive properties as for the
free Schr\"odinger equation.
This is very much reminiscent of the dispersive effects proved
for {\em small times} in \cite{Dsf} for wave equations with variable
coefficients, and relies on the fact that ${\wt \Psi}\approx -p^2/2$ as
$\eps t \leq \th$ is small enough.

In order to do so,
we need to get further informations both on the phase
$\wt\Psi$ and the amplitude $G$.

Firstly,
the smooth amplitude
$G$ is defined in (\ref{bigg}). It clearly is compactly supported
in $p$ and $\sqrt{\eps} q$.
Also, the gaussian $\exp(-x^2/2+i\G(\eps t,\sqrt{\eps}q,p) y \cdot y/2)$
belongs to the Schwartz space ${\cal S}\(\R^{2d}\)$ in the variables
$x$ and $y$ (recall indeed that ${\rm Im}\, \G(\eps t)>0$,
and $\eps t$ belongs to a compact set),
uniformly in the compactly supported parameters
$\eps t$, $\sqrt{\eps}q$, and $p$.
From this it follows that the amplitude
$G(q,p,\eta,\eps t,\sqrt{\eps}q,\sqrt{\eps}\eta)$
belongs to the Schwartz space
${\cal S}\(\R^{2d}\)$ in the first and third
variables $q$ and $\eta$, it is $C^\infty_c(\R^d)$
in the second variable $p$,
and these informations are uniform with
respect to the compactly supported
parameters $\eps t$, $\sqrt{\eps} q$, together with the
(non-compact) parameter $\sqrt{\eps}\eta$.

Secondly,
the smooth phase $\wt\Psi$
depends upon the small parameter $\eps t \in [0,2\th]$,
together with the two position/velocity variables
$\sqrt{\eps}q$ and $p$. All of them belong to a compact set.
It also depends upon the variable $\sqrt{\eps}\eta$, which is not
in a compact set.
On the more,
we have the easy first order expansion in the (small) parameter
$\eps t\leq 2 \th$,
\begin{align*}
&
\wt{\Psi}(p,\eps t,\sqrt{\eps}q,\sqrt{\eps}\eta)
=
\\
&
\qquad
-\frac{p^2}{2}+n^2(\sqrt{\eps}q)
-\sqrt{\eps} q \cdot \nabla_x n^2\(\sqrt{\eps}q\)
-\sqrt{\eps}\eta\cdot (p+O(\th))
+
O\(\th^2\)
.
\end{align*}
Here the remainder terms $O(\th)$, $O(\th^2)$, only depend upon
the compactly supported
parameters $\eps t\leq 2 \th$ and $p$, $\sqrt{\eps} q$ (they do not depend
upon $\sqrt{\eps} \eta$), and they are uniform with respect to these variables.
Hence, the stationary points of the phase (in the $p$ variable)
are given by
\begin{align}
\label{stpt}
&&
-p-\sqrt{\eps}\eta (1+O(\th))
+
O\(\th^2\)=0 .
\end{align}
Finally, there remains to observe
that the Hessian of the phase in $p$ is
\begin{align}
\label{hesspsi}
&&
\frac{D^2 \wt\Psi}{D p^2}
=
-{\rm Id}+O(\th) .
\end{align}

Upon taking $\th$ small enough, all these informations allow us to
make use of
the standard stationary phase estimate in $p$.
More precisely, we write,
\begin{align}
\label{fepsrev}
\non
&
f_\eps(t)
=
\int_{\R^{2d}}
\frac{dq d\eta}{\<q\>^{2d} \, \<\eta\>^{2d}}
\int_{\R^d} dp \,
\exp\(i \; t \; {\wt \Psi}(p,\eps t, \sqrt{\eps} q, \sqrt{\eps}\eta)\)
\\
&
\qquad\qquad\qquad
\<q\>^{2d} \, \<\eta\>^{2d}
\,
G(q,p,\eta,\eps t,\sqrt{\eps}q,\sqrt{\eps}\eta)
.
\end{align}
For each given values of $q$ and $\eta$,
we analyze the integral over $p$ in 
(\ref{fepsrev}). If $\sqrt{\eps} \eta$ is outside some compact
set around the support
of $G$ in $p$,
integrations by parts in
$p$ together with the information (\ref{stpt}), allow to prove
that the integral over $p$
in (\ref{fepsrev})
is bounded, for any integer $N$, by $C_{N, \th} t^{-N}$ for some
$C_{N, \th}>0$ independent of
$q$ and $\eta$. Hence the corresponding contribution to $f_\eps$
is bounded by
$C_{N, \th} t^{-N}$ as well.
Now, for $\sqrt{\eps}\eta$ in some compact set around the support
of $G$ in $p$,
we may use the information
(\ref{hesspsi}): this, together with the stationary phase Theorem with
the parameters
$\eps t$, $\sqrt{\eps} q$, $\sqrt{\eps} \eta$ in a compact set,
establishes that the 
integral over $p$ in (\ref{fepsrev}) is bounded by $C_\th t^{-d/2}$
for some $C_\th>0$, and
$C_\th$ turns out to be independent of
$q$ and $\eta$. Hence the corresponding contribuition to $f_\eps$
in (\ref{fepsrev}) is bounded
by $C t^{-d/2}$ as well.

All this gives the claimed estimate
\begin{align*}
&&
|f_\eps(t)|
\leq C_\th t^{-d/2} .
\end{align*}
The proof of proposition \ref{dur} is complete.


\section{Conclusion: Proof of the main Theorem}
\label{ccl}

\setcounter{equation}{0}


We want to prove
the convergence
$$
\<w^\eps,\phi\> \longrightarrow \<w^{\rm out},\phi\> ,
$$
when the source $S$ and the test function $\phi$ are Schwartz class.
Therefore, one needs to prove
\begin{align*}
&&
\frac{i}{\eps}
\int_0^{+\infty}
\e^{-\a_\eps t} \;
\<
U_\eps(t)
S_\eps
,
\phi_\eps
\>
\; dt
\rgt
\<w^{\rm out},\phi\> \; \text{ as } \eps \rgt 0 .
\end{align*}
Proposition \ref{ouane}
asserts
\begin{align*}
&
\frac{i}{\eps}
\int_0^{2 \T \eps}
\chi\(\frac{t}{\T \eps}\)
\
\e^{-\a_\eps t} \,
\<
U_\eps(t)
S_\eps
,
\phi_\eps\>
\, dt =
\<w^{\rm out},\phi\>
&
\\
&
\qquad\qquad\qquad\qquad\qquad
+O_{\T}(\eps^0)+O\(\frac{1}{T_0^{d/2-1}}\) ,
&
\end{align*}
where the notation $O(\eps^0)$ denotes a term going to zero with $\eps$,
and $O_{\T}(\eps^0)$ emphasizes the fact that the convergence depends
a priori on the value of
$\T$.

\noi
On the other hand
Proposition \ref{2} asserts
\beas
&&
\frac{1}{\eps}
\int_{\T \eps}^{+\infty}
\(1-\chi\)\(\frac{t}{\T \eps}\)
\;
\e^{-\a_\eps t} \;
\Big<
U_\eps(t)
\(1-\chi_\de\)\(H_\eps\) S_\eps
,
\phi_\eps\Big>
\; dt
\\
&&
\qquad\qquad\qquad\qquad
= O\(\frac{1}{\T}\)+O(\eps^0) .
\eeas
Now, for very large times and almost zero energies,
Proposition \ref{wawang} shows, for $\de$ small enough, and any $\k$,
\beas
&&
\frac{1}{\eps}
\int_{\eps^{-\k}}^{+\infty}
\e^{-\a_\eps t} \;
\Big<
U_\eps(t)
\chi_\de\(H_\eps\) S_\eps
,
\phi_\eps
\Big>
\; dt =O_{\k,\de}(\eps) .
\eeas
As for large times and almost zero energies, Proposition \ref{brob}
shows that, for $\de$ small enough,
$\k$ small enough, and $\tau$ large enough,
\beas
&&
\frac{1}{\eps}
\int_{\tau}^{\eps^{-\k}}
\e^{-\a_\eps t} \;
\Big<
U_\eps(t)
\chi_\de\(H_\eps\) S_\eps
,
\phi_\eps
\Big>
\; dt = O_{\k,\de}(\eps)
\eeas
Finally, for moderate times and almost zero energies, one has
the following two informations. First, for $\th$ small enough, and uniformly
in $\eps$, we have
\beas
&&
\frac{1}{\eps}
\int_{\T \eps}^{2\th}
\(1-\chi\)\(\frac{t}{\T \eps}\)
\chi\(\frac{t}{\th}\) \;
\e^{-\a_\eps t} \;
\Big<
U_\eps(t)
\chi_\de\(H_\eps\) S_\eps
,
\phi_\eps
\Big>
\; dt
\\
&&
\qquad\qquad\qquad\qquad
=O_\th\(\frac{1}{\T^{d/2-1}}\) .
\eeas
Second, for any fixed value of $\th>0$, and $\tau$,
\beas
&&
\frac{1}{\eps}
\int_{\th}^{\tau}
\(1-\chi\)\(\frac{t}{\T \eps}\)
\;
\e^{-\a_\eps t} \;
\Big<
U_\eps(t)
\chi_\de\(H_\eps\) S_\eps
,
\phi_\eps
\Big>
\; dt
\\
&&
\qquad\qquad\qquad
=
O_{\th,\tau,\de}\(\eps^0\) .
\eeas
All these informations show our main Theorem, upon conveniently
choosing the cutoff parameters
$\th$, $\T$, $\tau$ (in time), $\de$ (in energy), and the exponent
$\k$ (in time).
This ends our proof.


\section{Examples and counterexamples}
\label{expl}

\setcounter{equation}{0}



\subsection{The harmonic oscillator}


Given an appropriate potential $V(x)$, and
defining the semi-classical Schr\"odinger operator
$$
H_\eps=-\frac{\eps^2}{2} \D_x+V(x) ,
$$
our main Theorem proves
\bea
\label{cc}
\non
&&
\frac{1}{\eps}
\int_0^{+\infty}
\e^{-\a_\eps t}
\Bigg\<
\exp\( -i \frac{t}{\eps} H_\eps \)
S_\eps
,
\phi_\eps
\Bigg\> \, dt
\mathop{\longrightarrow}\limits_{\eps\rgt 0}
\\
&&
\qquad\qquad
\int_0^{+\infty}
\Big\<
\exp\( -i t \[-\D_x/2+V(0)\] \)
S
,
\phi
\Big\>
\, dt
.
\eea
Though we used in many places that our analysis requires
a potential of the form
$$
V(x)=
-n^2(x)=-n_\infty^2+O(\<x\>^{-\rho}) ,
$$
it seems interesting to investigate the validity of (\ref{cc})
when the potential is harmonic
\begin{align}
&&
V(x)=V(0)+\mathop{\sum}_{j=1}^d
\frac{\om_j^2}{2} x_j^2 ,
\end{align}
for some frequencies $\om_j \in \R$, and a given value $V(0)<0$.
Such a potential
does not enter our analysis since it is confining.
However, it is easily proved that for {\em pairwise
rationally independent} values of
the frequencies $\om_j$, the
transversality assumption {\bf (H)} page \pageref{HH} 
is true for this potential,
whereas in the extreme case where all $\om_j$'s are equal, this assumption
fails.
On the other hand,
one may use the Mehler formula \cite{Ho} (see \cite{C}
for a use of these formulae in the nonlinear context)
to compute the propagator
\begin{align}
\label{mehler}
&
\exp\(-i \frac{t}{\eps}\[-\eps^2 \D_x/2 + \sum_{j=1}^d \om_j^2 x_j^2/2 \]\)
=&
\\
&
\nonumber
\mathop{\prod}_{j=1}^d
\(
\frac{\om_j}{2 i \pi \eps \sin(\om_j t)}
\)^{1/2}
\exp\(\frac{i \om_j}{2 \eps \sin(\om_j t)}
\[(x_j^2+y_j^2)\cos(\om_j t)-2 x_j y_j\]\)
&
\end{align}
(Here we identified the propagator and its integral kernel).

Surprinsingly enough, using the Mehler formula to compute the limit
on the left-hand-side of (\ref{cc}), we may prove that for
{\em rationally independent} $\om_j$'s, the convergence result (\ref{cc})
is {\em locally true} in this case, for dimensions
$d\geq 4$, i.e (\ref{cc}) is true with the upper bounded $+\infty$ replaced 
by $T$, for any value of $T>0$.

We do not give the easy computations leading to this result.
The idea is the following: at each time $k\pi/\om_j$ ($k\in \Z$),
the trajectory of the harmonic oscillator shows periodicity
in the direction $j$. However, due to rational independence, at times
$k\pi/\om_j$, the trajectory does not show periodicity in any of the
$d-1$ other directions. Hence one gets enough local dispersion from
these directions to show that the corresponding contribution
to the time integral on the left-hand-side of (\ref{cc}) is roughly
$$
O\(
\int_{(-1+k\pi/\om_j)/\eps}^{(1+k\pi/\om_j)/\eps}
t^{-(d-1)/2}\, dt
\)
=
O\(
\eps^{(d-1)/2-1}
\)
\rgt 0 ,
$$
as long as $d-1>2$, i.e. $d\geq 4$.

Needless to say, in the extreme case where all $\om_j$'s are equal,
the result in (\ref{cc}) is {\em false}, even locally:
in this case, periodicity creates
a disastrous accumulation of energy at the origin ({\em all}
rays periodically hit the origin at times $k \pi/\om$, $k\in \Z$).

To our mind, this simple example indicates that our main Theorem
probably holds true for less stringent assumptions on the refraction index.
For instance, a uniform (in time) version of our transversality assumption
is probably enough to get the result (without assuming neither
decay at infinity of the refraction index,
nor assuming the non-trapping condition).


\subsection{Examples of flows satisfying the transversality condition}


We already observed that the harmonic oscillator with rationally independent
frequencies does satisfy the transversality assumption
{\bf (H)}. One actually has the value
$k=6d+1$ (see (\ref{93})) of the codimension in that case.

It is also easily verified
that the flow of a particle in a constant electric field,
i.e. the case of a potential
$$
V(x)=x_1 ,
$$
does satisfy {\bf (H)} as well, with $k=6d+1$.

Coupling the two flows, it is also verified that the potential
$$
V(x)=x_1+\sum_{j=1}^d \om_j^2 x_j^2/2 ,
$$
does satisfy {\bf (H)} as well, with $k=6d+1$.

Clearly, these examples are satisfactory, in that we may assume that
the potential has the above mentioned values {\em close} to the origin,
and we may truncate outside some neighbourhood of the origin so as to build
up a potential that satisfies the global assumptions we met in our main
Theorem.

\bigskip

\noi
{\bf Acknowledgements:}

{\small
\noi
The author wishes to thank B. Perthame and F. Nier for numerous
discussions on the subject. Also, he wishes to thank D. Robert for
his help concerning the wave-packet approach, and for giving him
the written notes \cite{Ro2}, as well as X.P. Wang for discussions and for
pointing out reference \cite{Wa}. Finally, he wishes to thank
R. Carles for a careful
reading of a first version of this manuscript.

This work has been partially supported by the ``ACI Jeunes Chercheurs -
M\'ethodes
haute fr\'equence pour les Equations diff\'erentielles ordinaires,
et aux d\'eriv\'ees partielles. Applications'',
by the GDR ``Amplitude Equations and 
Qualitative Properties'' (GDR CNRS 2103 : EAPQ) and the European Program 
'Improving the Human Potential' in the framework of the 'HYKE' network
HPRN-CT-2002-00282.
}


\begin{thebibliography}{99}
\bibliographystyle{alpha}



\bibitem[Ag]{Ag}
S. Agmon, {\it Spectral properties of Schr\"odinger
operators and scattering theory}, Ann. Scuola Norm. Sup. Pisa, Vol. 2,
N. 4, pp. 151-218 (1975).
 
\bibitem[AH]{AH}
S. Agmon, L. H{\"o}rmander,
{\it Asymptotic properties
of solutions of differential equations with simple characteristics},
J. Anal. Math. Vol. 30, pp. 1-37 (1976).

\bibitem[BCKP]{BCKP}
J.D. Benamou, F. Castella,
Th. Katsaounis, B. Perthame,
{\it
High frequency limit in the Helmholtz
equation},
Rev. Mat. Iberoamericana, Vol. 18, N. 1, pp. 187-209 (2002),
and,
S\'eminaire E.D.P., \'Ecole Polytechnique,
expos\'e N. V, 27 pp., 1999-2000.


\bibitem[BR]{BR}
A. Bouzouina, D. Robert,
{\it
Uniform semiclassical estimates for the propagation of quantum observables},
Duke Math. J., Vol. 111, no. 2, pp. 223-252 (2002).

\bibitem[Bu]{Bu}
N. Burq,
{\it
Semi-classical estimates for the resolvent in nontrapping geometries},
Int. Math. Res. Not., no. 5, pp. 221-241 (2002).

\bibitem[Bt]{Bt}
J. Butler,
{\it Global H Fourier Integral operators with complex valued phase functions},
preprint Bologna (2001).

\bibitem[C]{C}
R. Carles,
{\it
Remarks on nonlinear Schrödinger equations
with harmonic potential},
Ann. Henri Poincar\'e, Vol. 3, N. 4, pp 757-772 (2002).
 
\bibitem[CPR]{CPR}
F. Castella, B. Perthame, O. Runborg,
{\it High frequency limit of the Helmholtz equation II:
source on a general smooth manifold},
Comm. P.D.E., Vol. 27, N. 3-4, pp. 607-651 (2002).

\bibitem[CRu]{CR}
F. Castella, O. Runborg,
In preparation.

\bibitem[CRo]{CRo}
M. Combescure, D. Robert,
{\it
Semiclassical spreading of quantum wave packets
and applications near unstable fixed points of the classical flow},
Asymptot. Anal., Vol. 14, no. 4, pp. 377-404 (1997). 

\bibitem[CRR]{CRR}
M. Combescure, J. Ralston, D. Robert,
{\it
A proof of the Gutzwiller semiclassical trace formula using
coherent states decomposition},
Comm. Math. Phys., Vol.  202, no. 2, pp. 463-480 (1999). 


\bibitem[DG]{DG}
J. Derezi\'nski,
C. G\'erard,
{\bf
Scattering theory of classical and quantum $N$-particle systems},
Texts and Monographs in Physics, Springer-Verlag, Berlin (1997).

\bibitem[DS]{DS}
M. Dimassi, J. Sj\"ostrand,
{\bf
Spectral Asymptotics in the Semiclassical Limit},
London Math. Soc. Lecture Notes series 268, Cambridge Universuty Press
(1999).


\bibitem[Dsf]{Dsf}
D. Dos Santos Ferreira,
{\it In\'egalit\'es de Carleman $L^p$
pour des indices critiques},
Thesis, Universit\'e de Rennes 1 (2002).

\bibitem[Fo]{Fo}
G.B. Folland,
{\bf
Harmonic analysis in phase space}, Annals of Mathematics Studies, Vol. 122.
Princeton University Press, Princeton, NJ (1989).

\bibitem[Fou]{Fou}
E. Fouassier,
Thesis, University Rennes 1, in preparation.

\bibitem[GM]{GM}
C. G\'erard, A. Martinez,
{\it
Principe d'absorption limite pour des op\'erateurs de Schr\"odinger
\`a longue port\'ee},
C. R. Acad. Sci. Paris S\'er. I Math., Vol. 306, no. 3, pp. 121-123 (1988).

\bibitem[HJ]{HJ}
G.A. Hagedorn, A. Joye,
{\it
Semiclassical dynamics with exponentially small error estimates},
Comm. Math. Phys., Vol. 207, pp. 439-465 (1999).

\bibitem[H1]{H1}
G.A. Hagedorn,
{\it
Semiclassical quantum mechanics III},
Ann. Phys., Vol. 135, pp. 58-70 (1981).

\bibitem[H2]{H2}
G.A. Hagedorn,
{\it
Semiclassical quantum mechanics IV},
Ann. IHP, Vol. 42, pp. 363-374 (1985).

\bibitem[HR]{HR}
B. Helffer, D. Robert,
{\it Caclcul fonctionnel par la transformation de Mellin et
op\'erateurs admissibles},
J. Funct. Anal., Vol. 53, N. 3, pp. 246-268 (1983).

\bibitem[He]{He}
K. Hepp,
{\it
The classical limit of quantum mechanical correlation functions},
Comm. Math. Phys., Vol. 35, pp. 265-277 (1974).

\bibitem[Ho]{Ho}
L. H\"ormander,
{\it
Symplectic classification of quadratic forms, and general Mehler formulas},
Math. Z. Vol. 219, N. 3, pp. 413-449 (1995).
 
\bibitem[J]{J}
T. Jecko,
{\it
From classical to semi-classical non-trapping behaviour:
a new proof},
Preprint University of Rennes 1 (2002).

\bibitem[Ma]{Ma}
A. Martinez,
{\bf An introduction to Semicassical and Microlocal Analysis},
Universitext, Springer-Verlag, New York (2002).

\bibitem[PV1]{PV1}
B. Perthame, L. Vega,
{\it
Morrey-Campanato estimates for Helmholtz Equation},
J. Funct. Anal., Vol. 164, N. 2, pp. 340-355 (1999).

\bibitem[PV2]{PV2}
B. Perthame, L. Vega,
{\it Sommerfeld radiation condition for Helmholtz equation
with variable index at infinity}, Preprint (2002).


\bibitem[Ro]{Ro}
D. Robert,
{\it
Remarks on asymptotic solutions for time dependent Schr\"o\-din\-ger
equations},
to appear.

\bibitem[Ro2]{Ro2}
D. Robert, 
notes on lectures given at the University of Nantes (1999-2000).

\bibitem[Rb]{Rb}
S. Robinson,
{\it
Semiclassical mechanics for time-dependent Wigner functions},
J. Math. Phys., Vol. 34, pp. 2150-2205 (1993).

\bibitem[Wa]{Wa}
X.P. Wang,
{\it
Time decay of scattering solutions and resolvent estimates
for Semiclassical Schr\"odinger operators},
J. Diff. Eq., Vol. 71, pp. 348-395 (1988).

\bibitem[WZ]{WZ}
X.P. Wang, P. Zhang,
{\it
High frequency limit of the Helmholtz equation
with variable refraction index},
Preprint (2004).

\end{thebibliography}
\end{document}